\documentclass[a4paper,11pt]{article}
\usepackage{amsmath}
\usepackage{amsfonts}
\usepackage{titlesec}
\usepackage{amssymb}
\usepackage{amsthm}

\usepackage{geometry}
\usepackage{color}

\newcommand{\uzero}{u_0}

\usepackage[utf8]{inputenc}

\usepackage{mathtools}

\DeclareMathOperator{\sech}{sech}

\usepackage{array}
\makeatletter
\newcommand{\thickhline}{%
    \noalign {\ifnum 0=`}\fi \hrule height 1pt
    \futurelet \reserved@a \@xhline
}
\newcolumntype{"}{@{\hskip\tabcolsep\vrule width 1pt\hskip\tabcolsep}}
\makeatother

\newcolumntype{?}{!{\vrule width 1pt}}

\newcommand{\ci}{\mathrm{i}} %% complex number "i"

\usepackage{framed}
\usepackage{graphicx, caption, subcaption}
\usepackage{hyperref}
\usepackage{dsfont}
\usepackage{float}

\newtheorem{theorem}{Theorem}[section]

\newtheorem{remark}[theorem]{Remark}

\usepackage{diagbox}

\renewcommand{\thefootnote}{\fnsymbol{footnote}}
\renewcommand{\thefootnote}{\arabic{footnote}}

\begin{document}

\begin{center}
{\LARGE Numerical comparison of mass-conservative schemes for the Gross-Pitaevskii equation\renewcommand{\thefootnote}{\fnsymbol{footnote}}\setcounter{footnote}{0}
 \hspace{-3pt}\footnote{The authors acknowledge the support by the Swedish Research Council (grant 2016-03339) and the G\"oran Gustafsson Foundation}}\\[2em]
\end{center}

\renewcommand{\thefootnote}{\fnsymbol{footnote}}
\renewcommand{\thefootnote}{\arabic{footnote}}

\begin{center}
{\large Patrick Henning and Johan W\"arneg{\aa}rd \footnote[1]{Department of Mathematics, KTH Royal Institute of Technology, SE-100 44 Stockholm, Sweden.}}\\[2em]
\end{center}

\begin{center}
{\large{\today}}
\end{center}

\begin{center}
\end{center}

\begin{abstract}
In this paper we present a numerical comparison of various mass-conservative discretizations for the time-dependent Gross-Pitaevskii equation. We have three main objectives. First, we want to clarify how purely mass-conservative methods perform compared to methods that are additionally energy-conservative or symplectic. Second, we shall compare the accuracy of energy-conservative and symplectic methods among each other. Third, we will investigate if a linearized energy-conserving method suffers from a loss of accuracy compared to an approach which requires to solve a full nonlinear problem in each time-step. In order to obtain a representative comparison, our numerical experiments cover different physically relevant test cases, such as traveling solitons, stationary multi-solitons, Bose-Einstein condensates in an optical lattice and vortex pattern in a rapidly rotating superfluid. We shall also consider a computationally severe test case involving a pseudo Mott insulator. Our space discretization is based on finite elements throughout the paper. We will also give special attention to long time behavior and possible coupling conditions between time-step sizes and mesh sizes. The main observation of this paper is that mass conservation alone will not lead to a competitive method in complex settings. Furthermore, energy-conserving and symplectic methods are both reliable and accurate, yet, the energy-conservative schemes achieve a visibly higher accuracy in our test cases. Finally, the scheme that performs best throughout our experiments is an energy-conserving relaxation scheme with linear time-stepping proposed by C. Besse (SINUM,42(3):934--952,2004).
\end{abstract}


\section{Introduction}
This article deals with numerical methods for the nonlinear Schr\"odinger equation (NLSE) with cubic nonlinearity, which reads
		\begin{gather} \label{NLSE}
		\ci \hspace{1pt} \partial_t u = -
		\Delta u + Vu+\beta|u|^2u. 
		\end{gather}
		The solution $u(x,t)$ is a complexed valued wave function, $\beta \in \mathbb{R}$ is a given constant typically describing particle interactions and $V(x):\mathbb{R}^d \supset \Omega \mapsto \mathbb{R}$ represents a potential. The problem is completed by an initial value $u(x,0)=\uzero(x)$. The equation has applications in several fields such as optics \cite{Optics,GravityWaves}, fluid dynamics, in particular deep water waves \cite{FluidReview,DeepWater}, and quantum physics. Notably in quantum mechanics, it describes the dynamics of a Bose-Einstein condensate \cite{Dynamics}. The equation is often referred to as the Gross-Pitaevskii equation.
		
When selecting a suitable discretization for problem \eqref{NLSE}, one needs to be aware that the equation possesses several time invariant characteristics. Two of the most important are the mass $M$ and the energy $E$,
\begin{align}
&M[u] = \int_\Omega |u(x,t)|^2 \hspace{2pt} dx,  \label{Mass} \\
& E[u] = \int_\Omega 
|\nabla u(x,t)|^2 +V(x)|u(x,t)|^2 +\frac{\beta}{2}|u(x,t)|^{4} \hspace{2pt}dx \label{Energy}
\end{align}
which are preserved for all times. Here we stress that the selection of a numerical method that preserves mass and energy on the discrete scale is subject to the choice of a suitable time discretization, but it is typically not affected by the choice of the spatial discretization.

Another important characteristic of the problem is the Hamiltonian structure of the PDE. More precisely, by introducing $q:=\text{Re}(u)$ and $p:=\text{Im}(u)$ we observe that the NLSE can be interpreted as an infinite-dimensional Hamiltonian system of the form 
$$
	\partial_t q =  \frac{\delta H(p,q)}{\delta p} \qquad \mbox{and} \qquad
	\partial_t p =  -\frac{\delta H(p,q)}{\delta q},
$$
where $H(p,q)=\frac{1}{2}E[u]$ is the Hamiltonian defined through the energy functional. Discretizing the problem in space while staying continuous with respect to the time variable leads to a finite dimensional Hamiltonian system. For instance, if $p_i $ and $q_i$ represent finite difference values then the Hamiltonian system is of the form $\partial_t q_i = \partial H_h /\partial p_i$ and $\partial_t p_i = -\partial H_h / \partial q_i$ , where $H_h$ represents a corresponding discrete Hamiltonian. A second order finite difference method in one dimension would take the form
$$H_h = \frac{1}{2}\sum_i
 \left(\frac{q_i-q_{i-1}}{\Delta x} \right)^2+ 
 \left(\frac{p_i-q_{i-1}}{\Delta x} \right)^2+\frac{\beta}{2}(p_i^2+q_i^2)^2+V_i(q_i^2+p_i^2).$$
The following symplectic form is then time invariant: 
\begin{equation} \label{Symplecticity}
 \omega= \sum_i dp_i\wedge dq_i  . \end{equation}
It is well-known that symplectic time integrators only allow for small oscillations in the discrete energy, which is why they can be considered  almost energy conservative.
 
A central focus when solving the NLSE numerically has been to mimic time invariants, where the above mentioned invariants of mass, energy and symplecticity are considered to be the most important ones. An early comparative study was made by Sanz-Serna and Verwer \cite{Sanz-Serna} who compared closely related temporal discretizations based on finite differences in space. It was found that non mass-conservative methods were prone to nonlinear blow-up and performed poorly.  Moreover it was found that, among their test cases, the best-performing method was symplectic and mass-conservative. 
A slight modification of this method leads to an energy-conservative Crank-Nicolson method as shown by Sanz-Serna \cite{Sanz-SernaNLCN}. However, the conservation of energy comes at the cost of losing symplecticity. In fact, for $\beta \not = 0 $ and $d>1$ it is impossible to discretely conserve all of the three aforementioned time invariants simultaneously, since energy conservation and symplecticity are mutually exclusive if the scheme is not exact (cf. \cite{ZHONG1988134}). 
Other numerical comparisons, most often of widely different approaches, have been done notably by Chang et al. \cite{ComparisonFD} and Bao et al. \cite{BaoNum}.  A recent overview of the most common methods was done by Antoine et al. \cite{ANTOINE20132621}. Here we mention that the comparisons in \cite{ComparisonFD} and \cite{ANTOINE20132621} are only for $1d$ problems with smooth solutions, which do not allow representative conclusions to problems in $2d$ and with reduced regularity. For instance, the method that was found to be most robust and efficient in \cite{ComparisonFD} (Linearized Crank-Nicolson, cf. Section \ref{subsection-LCN-FEM}) can be affected by energy blow-ups in $2d$ (cf. Section \ref{EXP2D}).
 
In the present paper we complement the previous comparisons by answering which of the time invariants \eqref{Energy}, \eqref{Symplecticity}, in addition to mass \eqref{Mass}, is most important to conserve and how efficiently this can be done. A further aim is to find competitive methods to use in a low regularity regime and for $d>1$. We will fix the space discretization as being of finite element type (FEM) in order to clearly analyze and isolate the role of the time integrator. Changing the space discretization over the experiments could distort the picture, which is why we will avoid this (with one illustrative exception in Section \ref{StationarySoliton}).

The consideration of low regularity regimes is motivated by the non-smooth potentials that arise in the context of disorder in quantum systems where one, for instance, encounters random potentials, see \cite{DisorderQuantum} and \cite{DisorderLocal}. Another context is the investigation of quantum phase transitions from superfluid BEC to bosonic Mott insulators.

We stress that alternatively to a finite element discretization in space, one could also use a spectral or a finite difference discretization. In particular in high regularity regimes with smooth potentials and regular domains, the exponential convergence of a spectral approach leads to methods that are computationally very efficient. On the downside, they typically perform poorly if the regularity of the exact solution drops. Due to its popularity, we include the 2nd order Strang splitting spectral method \cite{Spectral} in one of our test cases which highlights these differences in characteristic features. For an overview of spectral methods for the Gross-Pitaevskii equation we refer to Bao et al. \cite{Spectral}. As for finite difference schemes, they are typically less effected by regularity issues. However, their accuracy still reduces observably.
An early paper by Griffiths \cite{GRIFFITHS1984177} found that finite element methods gave more favorable results than finite differences. A final argument for using FEM in connection with rough potentials is that it can be naturally combined with adaptive mesh refinement strategies, which may be key to efficient methods in low regularity regimes. For instance, rough disorder potentials are closely related to the phenomenon of Anderson localization \cite{AHP18,And58,ArnDJMF16} (exponential localization of eigenmodes), where often only locally refined meshes are required. For these reasons we restrict our attention to methods compatible with a finite element formulation. Concerning the application of nonlinear Schr\"odinger equations, we deem exponential time integrators to have low compatibility with FEM discretizations as this requires a mass lumping approach to be efficiently implemented. Since lumping can introduce a severe perturbation of the properties that are to be conserved, we exclude such methods in our survey.

All methods that we consider are mass-conservative and implicit. The necessity of mass-conserving schemes to obtain reliable approximations was demonstrated in the extensive numerical studies in \cite{Sanz-Serna}. Here, we also refer to the numerical experiments in \cite{HeM17} where the Backward Euler method shows a devastating performance due to its lack of mass conservation. If the nonlinearity is treated implicitly, then a nonlinear system of equations must be solved in each time-step.  This is the case for the symplectic one-stage Gauss-Legendre Runge-Kutta method that Sans-Serna and Verwer \cite{Sanz-Serna} found to perform best and also for the energy conservative Crank-Nicolson method.

Albeit implicit methods, convergence rates were, until recently, always obtained under a (counterintuitive) time-step condition in terms of the mesh size. 
Recent developments include analysis of linear mass-conservative methods by Wang \cite{Wang} (2013), Besse \cite{Besse} (2004), and Zouraris \cite{zouraris_2001} (2001) whose methods we take into account in our comparison. Wang proposes an Adams-Bashforth like linearization of the Crank-Nicolson method and was first to prove optimal convergence rates in the $L^2$-norm without any condition on the time-step size. The method proposed by Besse conserves, in addition to mass, energy. The Two-Step method analyzed by Zouraris was proved to converge under moderate time-step restrictions. As mentioned, the linearized Crank-Nicolson method was found best performing in \cite{ComparisonFD} for a set of $1d$ problems. The Relaxation method figures in the overview article by Antoine et al. \cite{ANTOINE20132621}. The authors present errors for a numerical experiment in which it, and the Crank-Nicolson method perform best in terms of temporal accuracy. Other than this, the linear methods have not, due to their recentness, figured in comparative studies and form an interesting complement to the well established nonlinear methods. 

Using a finite element based spatial discretization, we test these five mass-conservative methods on problems of increasing difficulty in one and two dimensions. The test cases include solitons with known analytical solution such as the single traveling soliton and a newly derived stationary soliton \cite{Exact2010}. Problems of physical relevance involving optical lattices and angular momentum are also considered. We do not resort to the usage of an artificial source function often found in the literature, as it would distort the physical behavior and break the conservation of mass, energy and symplecticity.

The outline of the paper is as follows:
in Section \ref{METHODS} the aforementioned five methods are presented along with their properties and a brief evaluation with respect to both theory and implementation. A summarizing table concludes the section. In Section \ref{EXP1D} and Section \ref{EXP2D} we present numerical experiments in one and two dimensions. Finally, conclusions are presented in Section \ref{Conclusions}. 
  
\section{General problem formulation and methods} \label{METHODS}
In this section we will state the five numerical methods that are at the core of our numerical study. We present the methods in a FEM formulation in space along with their properties and computational complexity. For a proper mathematical description of the setup, we let $\Omega\subset \mathbb{R}^d$ denote a bounded computational domain and let $H^1_0(\Omega)=H^1_0(\Omega,\mathbb{C})$ denote the Sobolev space of weakly differentiable, square-integrable and compactly supported functions over the complex field. By $\overline{z}$ we denote complex conjugation of a complex number $z$ and by $\langle u,v\rangle = \int_{\Omega} u \overline{v}$ the $L^2$-inner product. With this, we consider the following initial value problem in weak form. Given a real valued bounded potential $V \in L^{\infty}(\Omega)$ and an initial value $\uzero \in H^1_0(\Omega)$, find $ u\in L^\infty([0,T],H^1_0(\Omega)) \text{ with } \partial_t u \in L^\infty([0,T],H^{-1}(\Omega) ) 
$ such that $u(\cdot,0) = \uzero$  and such that: 
\begin{equation} \label{NLSEOmega}
\langle \ci \partial_t u, u \rangle  =  
\langle \nabla u, \nabla v\rangle + \langle Vu,v\rangle+\beta\langle|u|^2u,v\rangle 
\end{equation}
for all $v\in H^1_0(\Omega)$ and a.e. $t\in (0,T]$.

\begin{remark}
If $\beta\ge 0$ and $V\ge 0$, then problem \eqref{NLSEOmega} has at least one solution. This solution is also unique if $d=1,2$. For $d=3$, a corresponding uniqueness result is still open. For a proof of these results, as well as a more general existence theory in the focusing regime $\beta<0$ we refer to the excellent textbook by Cazenave \cite{Cazenave}.
\end{remark}
  
In the following, we let $S_h$ denote the space of P1 Lagrange finite elements on a quasi-uniform simplicial mesh on $\Omega$ with mesh size $h$. We use standard notation for the time discretization, i.e. for final computational time $T=N\tau$ and $n = 0,\cdots,N$ we have $t_n = n\tau$ and $u^n(x) = u(x,n\tau)$. To facilitate reading we also define $u^{n+1/2} = (u^{n+1}+u^n)/2 $.

As for complexity, it is noted that nonlinear methods suffer from an additional drawback. Namely, if the nonlinear system of equations is solved through a Newton step, then it cannot be done in the complex field since there is a non-holomorphic mapping, in this case $z\mapsto|z|$. To make this aspect explicit, let $u_h = \sum_i \xi_i v^h_i(x)$ be the approximate solution. Then $\frac{\mbox{d}}{\mbox{d}\xi}|\sum_i \xi_i v^h_i|^2 $ cannot be calculated as a function of $\xi$ if $\xi \in \mathbb{C}$. The solution must therefore be split into real and imaginary part. This leads to a system of equations of size $(2m)^2$ where $m$ is the number of Lagrange nodes. Linear methods lead to systems of size $m^2$ (over $\mathbb{C}$) with $2m^2$ degrees of freedom.

The following methods are selected to complement the earlier numerical studies, in particular that by Sanz-Serna and Verwer \cite{Sanz-Serna}. Among the methods investigated in the aforementioned comparison paper, we only include the symplectic {\it Gauss-Legendre Runge-Kutta method} since it performed best in their study.

\subsection{Implicit Midpoint Method (IM-FEM) }

The first scheme in our comparative study belongs to a general class of Gauss--Legendre Runge--Kutta methods, where we shall only consider the lowest order case (implicit midpoint method with finite element space discretization - IM-FEM). For realizations of higher order we refer to \cite{KaM99,Sanz-Serna1988}. It has been shown by Sanz-Serna \cite{Sanz-Serna1988} that all Gauss--Legendre Runge--Kutta methods are symplectic, which is also why the following one stage realization is symplectic.

For a given discrete initial value $u_h^0\in S_h$ the IM-FEM approximation $u^{n}_{h,\tau} \in S_h$ to $u^{n}$ for $n=1,\cdots,N$ is given by the time discretization
	\begin{equation}\label{Symplectic}
		\ci \big\langle \frac{u^{n+1}_{h,\tau}-u^n_{h,\tau}}{\tau},v\big\rangle =  
		\big\langle\nabla u^{n+1/2}_{h,\tau},\nabla v \big\rangle  + \big\langle Vu^{n+1/2}_{h,\tau},v\big\rangle + \beta \big\langle |u^{n+1/2}_{h,\tau}|^2u^{n+1/2}_{h,\tau},v\big\rangle,
	\end{equation}
for all $v\in S_h$.
Besides being symplectic, the method is also mass-conservative (i.e. $\| u^n_{h,\tau} \|_{L^2(\Omega)} = \| u^0_h \|_{L^2(\Omega)}$ for all $n\ge 0$) as is easily seen by testing in \eqref{Symplectic} with $v=u^{n+1}_{h,\tau} + u^n_{h,\tau}$ and taking the imaginary part of the equation.
We stress again that this method showed the best performance in the numerical comparison by Sanz-Serna and Verwer \cite{Sanz-Serna}.
In order to solve the nonlinear algebraic equation given by \eqref{Symplectic}, an inner loop of Newton iterations can be applied. 

A first proof of convergence in the fully discrete finite element setting was given in 1984 \cite{Verwer1984} for $d=1$. Later, for arbitrary space dimension, optimal $L^2$-convergence rates of order $\mathcal{O}(\tau^2+h^{2})$ were proven by Akrivis et al. \cite{Akrivis1991} for $V=0$ and under the time-step condition $\tau = \mathcal{O}(h^{d/4})$. Optimal $L^{\infty}$- and $H^1$-error estimates for $u^N-u_{h,\tau}^N$ were obtained by Tourigny \cite{Tou91}, who recovers again the constraint $\tau=\mathcal{O}(h^{d/4})$ to ensure convergence. In \cite{HeM17} it was shown that for optimal $L^2$- and $H^1$-error estimates, the constraint can be relaxed to $\tau=\mathcal{O}(|\ln{h}|^{-(1+\varepsilon)/4})$ for $d=2$ and to $\tau=\mathcal{O}(h^{(1+\varepsilon)/4})$ for $d=3$. Here, $\varepsilon>0$ can be arbitrary close to zero. We note that the analysis in \cite{HeM17} also amounts for potentials $V\ge 0$, provided that they are smooth enough.

\subsection{Crank-Nicolson Method (CN-FEM)  }   
A slight modification of the previous method yields the energy-conservative Crank--Nicolson method (CN-FEM). The method has been of widespread use in the physics community and reads:
given $u_h^0\in S_h$ find $u^{n}_{h,\tau} \in S_h, \ n=1,\cdots,N$, such that for all $v\in S_h$:
	\begin{equation}\label{Crank}
		\ci \big\langle \frac{u^{n+1}_{h,\tau}-u^n_{h,\tau}}{\tau},v\big\rangle = 
		\big\langle\nabla u^{n+1/2}_{h,\tau},\nabla v \big\rangle + \big\langle Vu^{n+1/2}_{h,\tau},v\big\rangle + \beta \big\langle \frac{|u^{n+1}_{h,\tau}|^2+|u^{n}_{h,\tau}|^2}{2}u^{n+1/2}_{h,\tau},v\big\rangle.
	\end{equation}
Note that the essential difference between the methods \eqref{Symplectic} and \eqref{Crank} is that the IM-FEM is based on an average of the wave-functions in the nonlinearity, whereas the CN-FEM is based on an average of the densities.

\begin{remark}
A generalization of the IM-FEM to a larger class of nonlinearities is straightforward. To generalize the CN-FEM, antiderivatives of the nonlinearity are needed. For instance, a generalization of the Crank-Nicolson method for nonlinearities of the form $\gamma(|u|^2)u$ for some smooth and integrable $\gamma$ reads as follows in the fully discrete case:
	\begin{eqnarray*}
	\lefteqn{\ci \big\langle \frac{u^{n+1}_{h,\tau}-u^n_{h,\tau}}{\tau},v\big\rangle} \\
	&=& 
	\big\langle\nabla u^{n+1/2}_{h,\tau},\nabla v \big\rangle + \big\langle Vu^{n+1/2}_{h,\tau},v\big\rangle + \left\langle \frac{\Gamma(|u^{n+1}_{h,\tau}|^2)-\Gamma(|u^{n}_{h,\tau}|^2)}{|u^{n+1}_{h,\tau}|^2-|u^n_{h,\tau}|^2}u^{n+1/2}_{h,\tau},v\right\rangle
	\end{eqnarray*}
where \[\Gamma(\rho) := \int_0^\rho \gamma(t) \hspace{2pt} dt. \] 
\end{remark}

The method is mass and energy conservative. A numerical analysis of this method for $V=0$ was first done by Sanz-Serna in one space dimension \cite{Sanz-SernaNLCN}, in two and three dimensions by Akrivis et al. \cite{Akrivis1991}. In the latter work, optimal $L^2$-convergence rates, i.e. $\mathcal{O}(\tau^2+h^{2})$, were obtained under the time-step condition $\tau=\mathcal{O}(h^{d/4})$. Recently it was shown that such coupling conditions are in fact not required. More precisely, in \cite{NonlinearCN} optimal convergence in $L^2$ is proved without any time-step conditions for $d= 1,2,3$, for sufficiently smooth $V$ and a general class of nonlinearities including the cubic case $\gamma : z \mapsto |z|^2$. Furthermore for general potentials $V\in L^\infty(\Omega)$ the authors of \cite{NonlinearCN} prove suboptimal convergence rates of order $\mathcal{O}(h^{(d+\alpha)/2}+\tau)$ for some $\alpha>0$ under the mild condition $h^{4-d-\alpha}\leq C\tau^2$. We note that a proof of convergence in $H^1$ is still open in the literature. 

Practically, the CN-FEM in \eqref{Crank} requires the solving of a nonlinear system of equations in each time-step, a suitable Newton step was proposed and analyzed by Akrivis et al.

\subsection{Relaxation Method (RE-FEM)}
A relaxation method was put forward by Besse \cite{Besse} and treats the nonlinearity explicitly. A remarkable property of this method is that it conserve mass and energy simultaneously, while being linear in each time-step. Introducing $\rho_\tau = |u_\tau|^2$, the relaxation method in semi-discrete form reads as follows for $n\ge0$. Find $u_\tau^{n+1} \in H^1_0(\Omega)$ such that
\begin{align} \label{Relax}
\ci \frac{u_\tau^{n+1}-u_\tau^n}{\tau} = - 
\Delta u_\tau^{n+1/2} +V u_\tau^{n+1/2}+\beta \rho_\tau^{n+1/2} u_\tau^{n+1/2},
\end{align}
where the density approximation $\rho_\tau^{n+1/2}$ is recursively defined for $n\ge 0$ by
\begin{align*}
		 \frac{\rho_\tau^{n+1/2}+\rho_\tau^{n-1/2}}{2} = |u_\tau^n|^2 \qquad \mbox{and} \qquad \rho_\tau^{-1/2} := |\uzero|^2.
\end{align*}
Here, $u_\tau^n$ denotes the obtained numerical approximation for $u(\cdot,t_n)$ and $u_\tau^0=\uzero$ is the original initial value. To improve accuracy, an alternative initialization is to solve for $\rho_{\tau}^{-1/2}$ and then use the above recursion. 
The following initial step would then be added:
\begin{align}\label{Initialization2}
	\ci \frac{u_\tau^{0}-u_\tau^{-1/2}}{\tau} & = - 
	\Delta (u_\tau^{0}+u_\tau^{-1/2}) +V (u_\tau^{0}+u_\tau^{-1/2})+\beta |u^0_\tau|^{2} (u_\tau^{0}+u_\tau^{-1/2}) \\
	 \rho^{-1/2}_\tau & = |u^{-1/2}_\tau|^2 . \nonumber
\end{align}
 For the time-discrete system with $V=0$ and an appropriate choice for $\rho^{-1/2}_{\tau}$, 
optimal convergence rates were proved under the assumption of sufficiently high regularity \cite{BDD18}. More specifically, if the initial value fulfills $u_0 \in H^{4+s}(\mathbb{R}^d)$ for $s>d/2$, if $\| \rho^{-1/2}_{\tau} - |u(-\frac{\tau}{2})|^2 \|_{H^{4+s}(\mathbb{R}^d)} = \mathcal{O}(\tau^2)$ and if $t\mapsto u(t)$ is a smooth map where $u(t) \in H^{4+s}(\mathbb{R^d})$ solves
	\begin{equation} 
		\begin{dcases}
		\ci \partial_t u = - 
		\Delta u+\beta |u|^2u \\
		u(x,0) = \uzero(x),
		 \end{dcases}
	\end{equation}
then the sequence of time-discrete solutions $(u_\tau,\rho_\tau)$ to \eqref{Relax} converges to the exact solution $(u,|u|^{2})$ in $L^\infty([0,T];(H^s(\mathbb{R}^d))^2) $ with order $\mathcal{O}(\tau^2)$.
	
In the fully discrete case, the Relaxation Finite Element Method (RE-FEM) reads as follows in variational form.
Given a suitable approximation $u^{0}_h \in S_h$ to the exact initial value $\uzero$, find $u^{n}_{h,\tau} \in S_h, \ n=1,\cdots,N$, such that for all $v\in S_h$:
\begin{gather} 
\begin{dcases*}
 \frac{\rho^{n+1/2}_{h,\tau}+\rho^{n-1/2}_{h,\tau}}{2} = |u^n_{h,\tau}|^2 \\
\ci \big\langle \frac{u^{n+1}_{h,\tau}-u^n_{h,\tau}}{\tau},v\big\rangle =  
\big\langle\nabla u^{n+1/2}_{h,\tau},\nabla v \big\rangle  + \big\langle Vu^{n+1/2}_{h,\tau},v\big\rangle + \beta \big\langle \rho^{n+1/2}_{h,\tau} u^{n+1/2}_{h,\tau},v\big\rangle.
\end{dcases*}
\end{gather}
 The following energy like quantity is exactly conserved:
\begin{equation}
	E_0 = \int_{\Omega}  
	|\nabla u^{0}_h |^2 + \frac{\beta}{2}|u^{0}_h |^4 + V \hspace{2pt} |u^{0}_h |^2  = \int_{\Omega}  
	|\nabla u^n_{h,\tau}|^2 + \frac{\beta}{2}\rho^{n+1/2}_{h,\tau} \rho^{n-1/2}_{h,\tau} + V|u^n_{h,\tau}|^2 
\end{equation}
Note that if $u^n_{h,\tau}$ is as assumed piecewise linear, then for energy conservation it is imperative that the discrete density, $\rho^n_{h,\tau}$, be accurately represented as a piecewise polynomial of degree two. This needs to be considered when implementing the method. 

We are not aware of any analytical results concerning the convergence of the fully discrete method for the nonlinear Schr\"odinger equation, but for the semilinear heat equation an analysis was recently provided in \cite{Zou18}.

\subsection{Linearized Crank-Nicolson Method (LCN-FEM)}
\label{subsection-LCN-FEM}
The following method, first proposed and analyzed by Wang \cite{Wang}, is based on replacing the nonlinear term by an Adams--Bashforth linearization. For that we define the short hand notation
\begin{align*}
\widehat{u}^n_{h,\tau} := \frac{1}{2}(3u^{n}_{h,\tau}-u^{n-1}_{h,\tau}) \qquad \mbox{for } n\ge 1
\end{align*}
and for $n=0$, we let $\widehat{u}^0_h \in S_h$ denote the solution to the intermediate first step
\begin{align*}
\ci \big\langle \frac{\widehat{u}^0_{h,\tau} -u^0_h }{\tau/2}\big\rangle = \big\langle \nabla \widehat{u}^0_{h,\tau},\nabla v\big\rangle + 
 \big\langle V\widehat{u}^0_{h,\tau},v\big\rangle
+
\beta \big\langle |u^0_h|^2\widehat{u}^0_{h,\tau},v\big\rangle \qquad \mbox{for all } v \in S_h.
\end{align*}
With this, the Linearized Crank-Nicolson Finite Element Method (LCN-FEM) reads: find $u^n_{h,\tau} \in S_h$ such that 
\begin{align*}
\ci \big\langle \frac{u^{n+1}_{h,\tau}-u^n_{h,\tau}}{\tau},v\big\rangle =  
\big\langle\nabla u^{n+1/2}_{h,\tau},\nabla v \big\rangle +\big\langle Vu^{n+1/2}_{h,\tau},v\big\rangle + \beta \big\langle |\widehat{u}^n_{h,\tau}|^2 u^{n+1/2}_{h,\tau},v\big\rangle
\end{align*}
for all $v \in S_h$. Applying this particular Adams--Bashforth linearization can be motived by elliptic regularity theory: apart from the space discretization, the problem that needs to be solved for each time-step is a linear elliptic problem that admits $H^2$-regularity. Exploiting this observation allows to straightforwardly derive uniform $L^{\infty}$-bounds for the arising semi-discrete approximations. This a priori control over the growth of the numerical approximations is the key to a convergence theory that does not require time-step constraints. At the cost of this elegant approach, energy-conservation and symplecticity are lost. We stress that Wang \cite{Wang} was the first to prove optimal $L^2$-convergence rates $\mathcal{O}(\tau^2+h^{2})$ without coupling conditions for the mesh size and the time-step size. The analysis was carried out under the assumption that $V=0$. The paper also provides a proof of optimal order convergence in $H^1$, i.e. convergence with the rate $\mathcal{O}(\tau^2+h)$.

\subsection{Linearized Two-step Method (Two-Step FEM)}
In \cite{zouraris_2001}, Zouraris proposed the a new type of linearized method to which we will refer as Two-Step FEM. In the initial double step, the method requires to first solve for an intermediate step solution $u^{1/2}_{h,\tau} \in S_h$ with
\begin{align*}
 \ci \big\langle \frac{u^{1/2}_{h,\tau}-u^{0}_h}{\tau/2},v\big\rangle =  
 \big\langle\nabla\big(\frac{u^{1/2}_{h,\tau}+u^{0}_h}{2}\big),\nabla v \big\rangle + \big\langle V\frac{u^{1/2}_{h,\tau}+u^{0}_h}{2},v\big\rangle + \beta 
 \big\langle |u^{0}_h|^2 \hspace{2pt}\frac{u^{1/2}_{h,\tau}+u^{0}_h}{2},v\big\rangle
\end{align*}
for all $v \in S_h$. After that, the first full step is taken by solving for $u^{1}_{h,\tau} \in S_h$ with
\begin{align*}
\ci \big\langle \frac{u^{1}_{h,\tau}-u^{0}_h}{\tau},v\big\rangle =  
\big\langle\nabla (\frac{u^{1}_{h,\tau} + u^{0}_h}{2}),\nabla v \big\rangle + \big\langle V\frac{u^{1}_{h,\tau} + u^{0}_h}{2},v\big\rangle + \beta 
\big\langle |u^{1/2}_{h,\tau}|^2  \hspace{2pt} \frac{u^{1}_{h,\tau} + u^{0}_h}{2},v\big\rangle
\end{align*}
for all $v \in S_h$. With this, the iterations of the two-step method for $n=1,\cdots, N-1$ read: find $u^{n}_{h,\tau} \in S_h$ such that
\begin{eqnarray*}
\lefteqn{\ci \big\langle \frac{u^{n+1}_{h,\tau}-u^{n-1}_{h,\tau}}{2\tau},v\big\rangle} \\
&=&  
\big\langle\nabla (\frac{u^{n+1}_{h,\tau}+u^{n-1}_{h,\tau}}{2}),\nabla v \big\rangle + \big\langle V\frac{u^{n+1}_{h,\tau}+u^{n-1}_{h,\tau}}{2},v\big\rangle + \beta
\big\langle |u^{n}_{h,\tau}|^2 \frac{u^{n+1}_{h,\tau}+u^{n-1}_{h,\tau}}{2} ,v\big\rangle.
\end{eqnarray*}
for all $v \in S_h$.

In \cite{zouraris_2001}, the method was analyzed for the case of adaptive meshes, though we restrict our discussion of the results to the quasi-uniform case as before.
For $V=0$, optimal convergence rates in the $L^2$- and $H^1$-norm of order $\mathcal{O}(\tau^2+h^{2})$ and $\mathcal{O}(\tau^2+h)$ respectively are presented under certain constraints for the time-step size. For $d=1$ there is no constraint at all and for higher dimensions the time-step size needs to be bounded in terms of the mesh size by fulfilling the following relations
\begin{align*}
\tau=\mathcal{O}(|\ln{h}|^{-1/3}) \quad \text{for } d =2 \qquad \mbox{and} \qquad 
\tau=\mathcal{O}(h^{1/3}) \quad \text{for } d =3.
\end{align*}
We note that no such coupling constraints became visible in our numerical experiments.
 
 \subsection{Summarizing table}

The following table summarizes the essential features of the considered methods, i.e. mass and energy conservation, symplecticity and linear/nonlinear time-stepping. We recall that all methods have the same optimal convergence rates in space and time, provided that the solution is sufficiently regular. As for the detailed convergence results including the required assumptions we refer to the discussion in the corresponding subsections above. We stress again that the nonlinear CN-FEM is the only method for which convergence for rough potentials is proved. We do not include this fact in the table below. However, we shall indicate in the table if optimal $L^{\infty}(L^2)$- or $L^{\infty}(H^1)$-estimates are available for a certain method (provided sufficient regularity) and if the proofs for these convergence rates required a formal coupling condition between mesh size and time-step size. By {\it conditional} we mean that a coupling condition was required and by {\it unconditional} that it was not required. The precise conditions can be found in the previous subsections. For the Besse relaxation scheme, error estimates are only available for the semi-discrete method \cite{BDD18}. As possible coupling conditions and rates in $h$ remain open, we do not include it in the table.

\begin{table}[h!]
	\centering\scriptsize
 \begin{tabular}{ l"  c c c c c l }

 	\backslashbox{Property \ \ \ }{Method}  & IM-FEM & CN-FEM & RE-FEM & LCN-FEM & Two-Step FEM   \\ \thickhline
 	Mass-Conservative & Yes & Yes & Yes & Yes & Yes \\ \hline
 	Energy-Conservative & No & Yes & Yes  & No & No \\ \hline
 	Symplectic & Yes & No & No & No & No \\ \hline
 	Linear time-steps& No & No & Yes & Yes & Yes \\ \hline
 	$L^2$-convergence rates & Optimal  & Optimal & -  & Optimal & Optimal \\ 
	  & Conditional & Unconditional &  & Unconditional & Conditional \\ \hline
 	$H^1$-convergence rates & Optimal & - & - & Optimal & Optimal \\ 
 	  & Conditional &  &  & Unconditional & Conditional
 	\\ \hline
 	\end{tabular} 
 
 \caption{Note that an optimal convergence rate for the error $e_h = u_{h,\tau}^N - u^N$ in $L^2$ is of order $\mathcal{O}(\tau^2+h^2)$ and an optimal convergence rate in $H^1$ is of order $\mathcal{O}(\tau^2+h)$. As everywhere in the paper, this refers to the case of piecewise linear FEM (simplicial Langrange elements) for the space discretization.}
\end{table}

We stress that the computational complexity of the linearized methods RE-FEM, LCN-FEM and Two-Step FEM is basically the same (for a fixed space and time resolution). In comparison, the computational complexity of the nonlinear methods IM-FEM and CN-FEM is roughly ten times higher in one dimension. An exemplary comparison of CPU times is given in Table \ref{CPU1D} in Section \ref{StationarySoliton} . As we will see later, the higher CPU times can be justified for some test cases by a significantly higher accuracy. Furthermore, we note that even though two of the methods require coupling conditions between $h$ and $\tau$ for a formal proof of convergence, we could not find any evidence for this in our experiments. It seems that convergence is obtained for any mesh/step size ratio as long as $\tau,h \rightarrow 0$.

\section{Experiments in 1D} \label{EXP1D}
\subsection{Single Soliton}
In the first numerical experiment we consider a single traveling soliton $u$ given by
\begin{align}
\label{signle-solition}
u(x,t) = \sqrt{2}e^{\ci (\frac{x}{2}+\frac{3t}{4})}\sech(x-t).
\end{align}
It is characterized as the solution to the following NLSE in full space	
\begin{equation*}
\begin{cases}
\ci \partial_t u &= - \partial_{xx} u - |u|^2u  \qquad \mbox{in } \mathbb{R}\times (0,T] \\
u(x,0) &= \sqrt{2}e^{\frac{\ci x}{2}}\sech(x) \qquad\hspace{2pt} \mbox{in } \mathbb{R}.
\end{cases}
\end{equation*}
This solution originates from an early paper by Shabat and Zakharov \cite{SingleSoliton}. For any fixed time $t$, the solution given by \eqref{signle-solition} decays exponentially as $|x|\rightarrow \infty$. We consider the computational domain $[-30, 70]$ with homogeneous Dirichlet boundary conditions and maximal time $T=10$. This test case describes a simple setup without potential and with a solution that is smooth in space and time.

We run all five methods stated in Section \ref{METHODS} for various mesh sizes and time-step sizes. The discrete initial value is the nodal interpolation of $u(x,0)$ in the finite element space. As we will see, all methods perform very well for the test case, where the Two-Step FEM performs best in terms of error size and CPU times. 

\begin{figure}[H]
	\vspace{-40pt}
		\centering

		\begin{subfigure}[t]{0.4\textwidth}
			\includegraphics[width=1.0\linewidth]{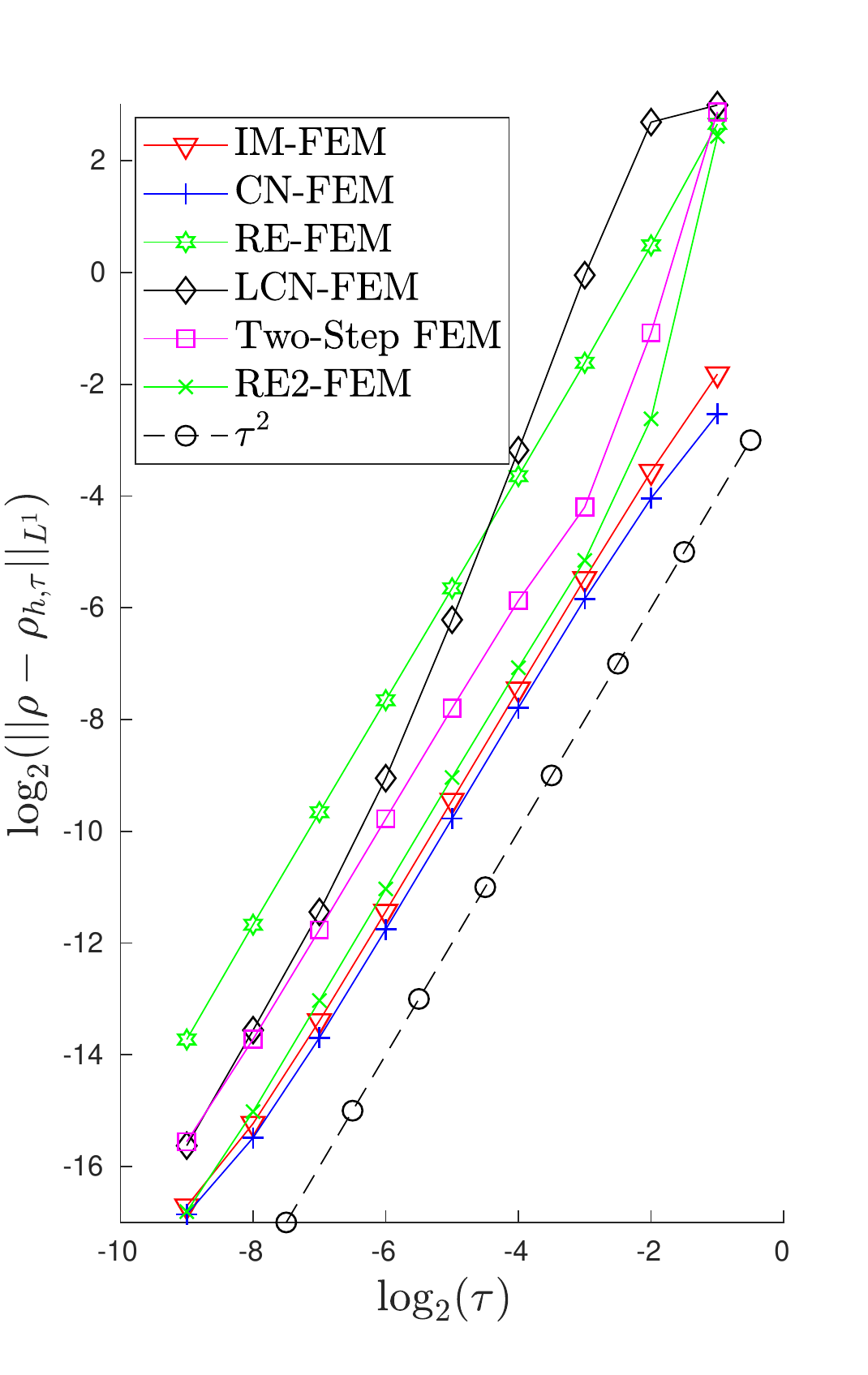}
			\caption{\scriptsize $L^1$-norm of error in density for the spatial mesh size $h=2^{-10}$ plotted versus $\tau$.}
			\label{Test1Rho}
		\end{subfigure}
		 \hspace{-5pt}
		~
		\begin{subfigure}[t]{0.4\textwidth}
			\includegraphics[width=1.0\linewidth]{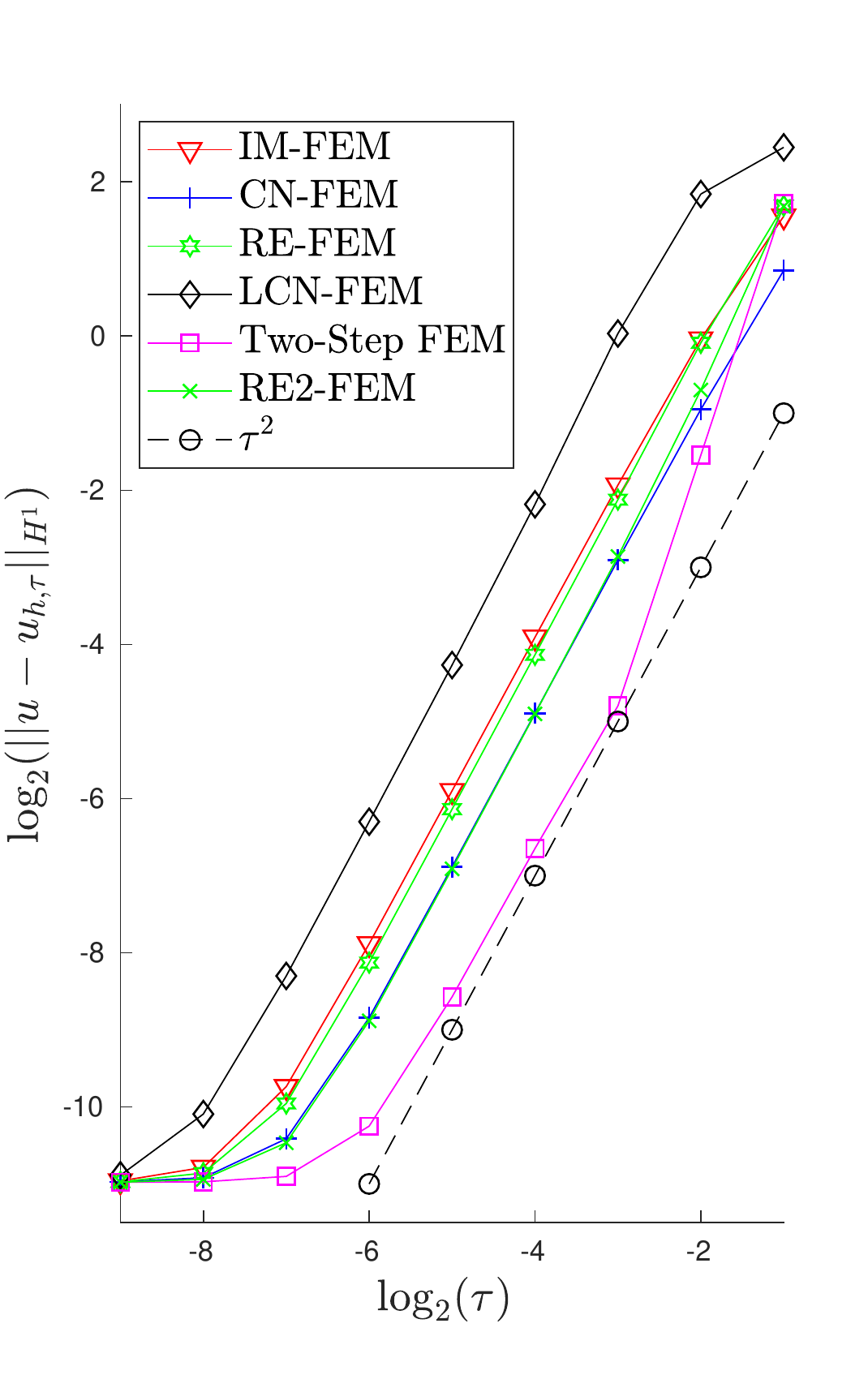}
			\caption{ \scriptsize $H^1$-norm of error for the spatial mesh size $h=2^{-10}$ plotted versus $\tau$.}
			\label{Test1H1}
		\end{subfigure}%
		 \hspace{-5pt}
		 \caption{$L^1$-errors of density and $H^1$-errors for all five methods for the single soliton test case.}
\end{figure}

To support this claim, a first set of results is depicted in Fig. \ref{Test1Rho}. Here, the $L^1$-error for the density $\rho = |u|^2$ at time $T=10$ is plotted versus the time-step size $\tau$. The nonlinear methods start with substantially smaller errors than the linear methods. Overall the nonlinear methods outperform the linear methods for coarse time steps, where the energy-conservative CN-FEM has only slightly smaller errors than the symplectic IM-FEM. For smaller time steps, the linear methods catch up in terms of performance due to their shorter CPU times (making them 10 times faster than the nonlinear methods). We observe that the LCN-FEM initially converges cubically in the $L^1$ norm of the density, it does however start with the highest errors. At $\tau = 2^{-7}$ it is on par with the Two-Step FEM and the cubic convergence flattens out to the expected quadratic convergence. 
Asymptotically, the Two-Step FEM and the LCN-FEM are on par, with roughly  four times larger errors than the nonlinear methods. The RE-FEM in turn produces about four times larger errors than the Two-Step FEM. The relatively large errors of the RE-FEM are due to the initialization $\rho^{-1/2} = |u_0|^2$, which causes a larger initial error. If instead $\rho^{-1/2}$ is solved for using \eqref{Initialization2}, then the initial error is a small perturbation of order $\mathcal{O}(\tau^2)$. In this case, the $L^1$ error of the density becomes on par with those of the nonlinear methods and the $H^1$-error almost exactly that of the CN-FEM. This is plotted in Figures \ref{Test1Rho} and \ref{Test1H1} as RE2-FEM.

In Fig. \ref{Test1H1} the error in the $H^1$-norm is plotted versus the time-step size. Again, the Two-Step FEM fares best. Next in accuracy comes the  energy conservative CN-FEM and only slightly less accurate is the energy conservative RE-FEM which is on is on par with the symplectic IM-FEM. The LCN-FEM consistently has about 6 times higher errors than the CN-FEM and roughly 32 times that of the Two-Step FEM.

It is noted in Fig. \ref{Test1H1} that for approximately $\tau = 2^{-9}$ the error in $H^1$-norm from the space discretization with $h=2^{-10}$ seems to become dominant in all methods, resulting in polluted convergence rates in $\tau$. 

\clearpage

	\begin{figure}[H]
		\centering

		\begin{subfigure}[t]{0.4\textwidth}
			\includegraphics[width=1.0\linewidth]{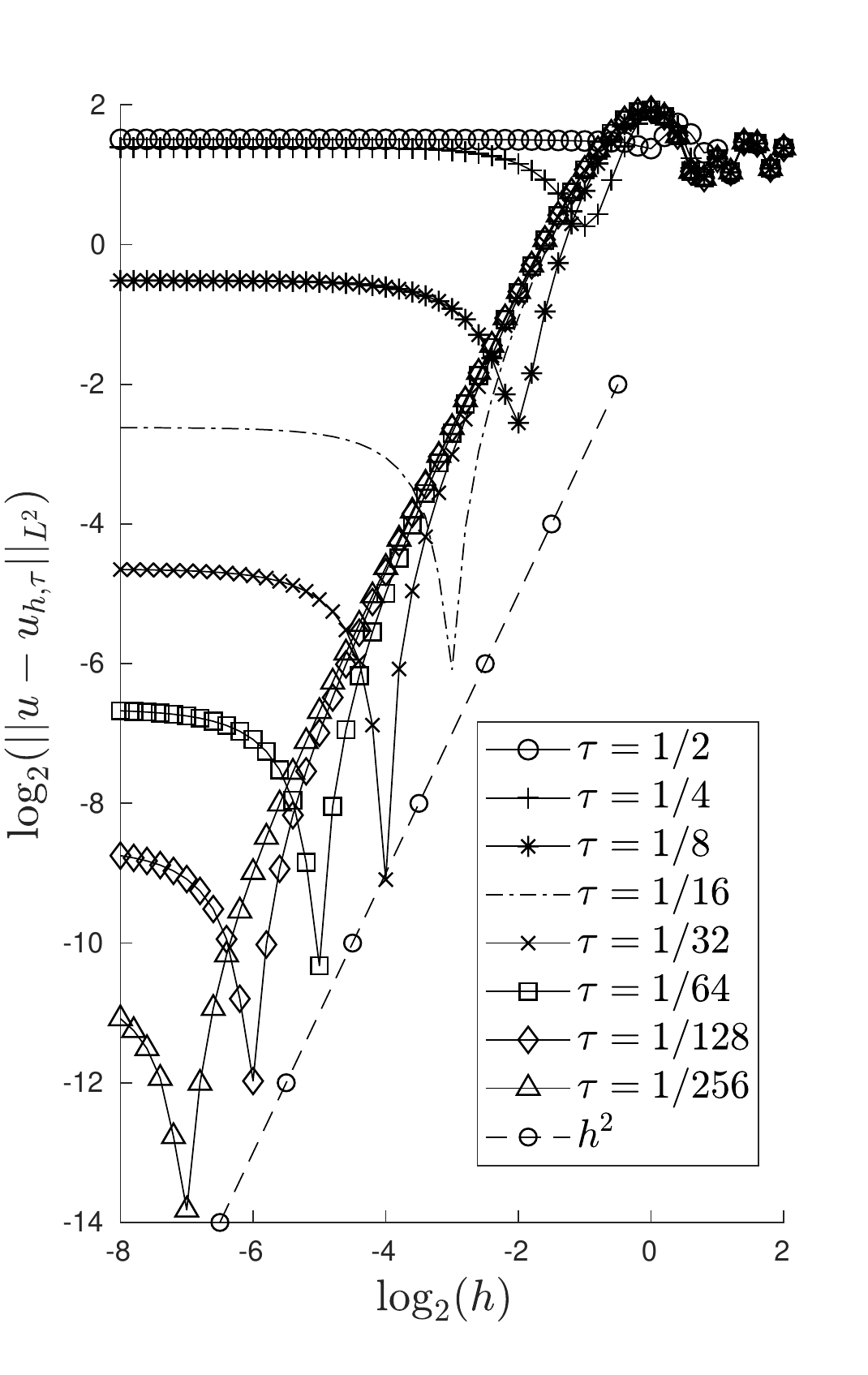} \vspace{-15pt}
			\caption{\scriptsize $L^2$-errors for the LCN-FEM.  }
						\label{Test1WangL2}
		\end{subfigure}
		~
		\begin{subfigure}[t]{0.4\textwidth}
			\includegraphics[width=1.0\linewidth]{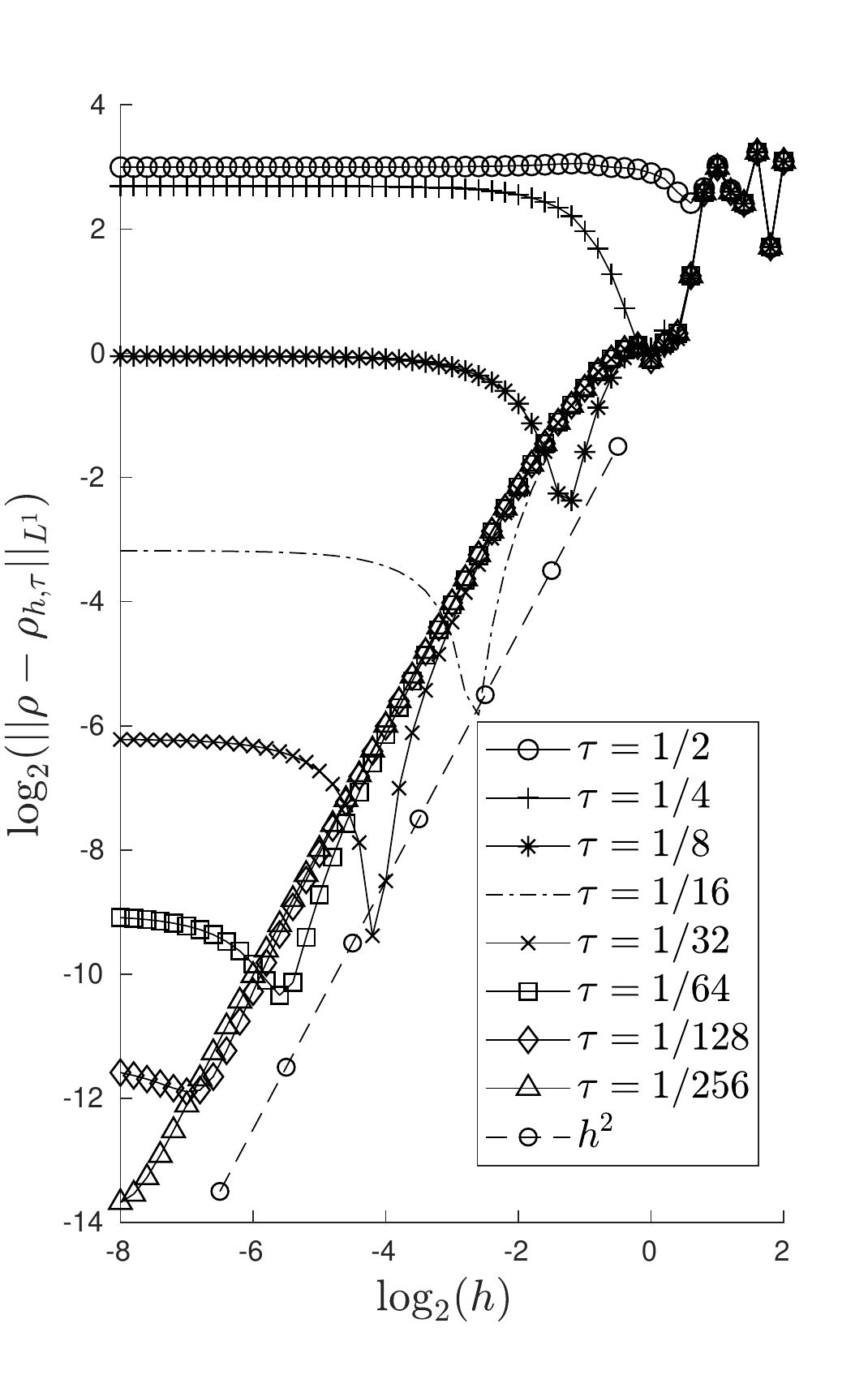}
			\caption{ \scriptsize $L^1$-errors of density for the LCN-FEM.}
			\label{WangRHO}
		\end{subfigure}%
	\caption{We observe remarkably small $L^2$-errors for the LCN-FEM when the spatial and temporal discretizations are coupled as $h/\tau =2$. This feature virtually disappears in the error plot of the density, implying that for this ratio the error in phase is minimal.}
\end{figure}

	\clearpage

An interesting observation is made for the LCN-FEM, which achieves remarkably low $L^2$-errors if the mesh size and the time-step size are coupled to the ratio $h/\tau = 2$. In this case, the method still shows a quadratic convergence, but the errors are significantly smaller. This is illustrated in Fig. \ref{Test1WangL2}. The ratio does not depend on final time and is not observed for the other methods. As seen in Fig. \ref{WangRHO} this feature virtually disappears in the error plot of the density, implying that for this ratio the error in phase is minimal. That is, a small change in the spatial discretization may result in considerably smaller  $L^2$-errors, but not in smaller density errors. Since the density barely changes, the phase must account for the smaller errors.	
	
	\subsection{Two Stationary Solitons} \label{StationarySoliton}
\begin{figure}[h]
    \centering
    \begin{subfigure}[b]{0.3\textwidth}
        \includegraphics[width=\textwidth]{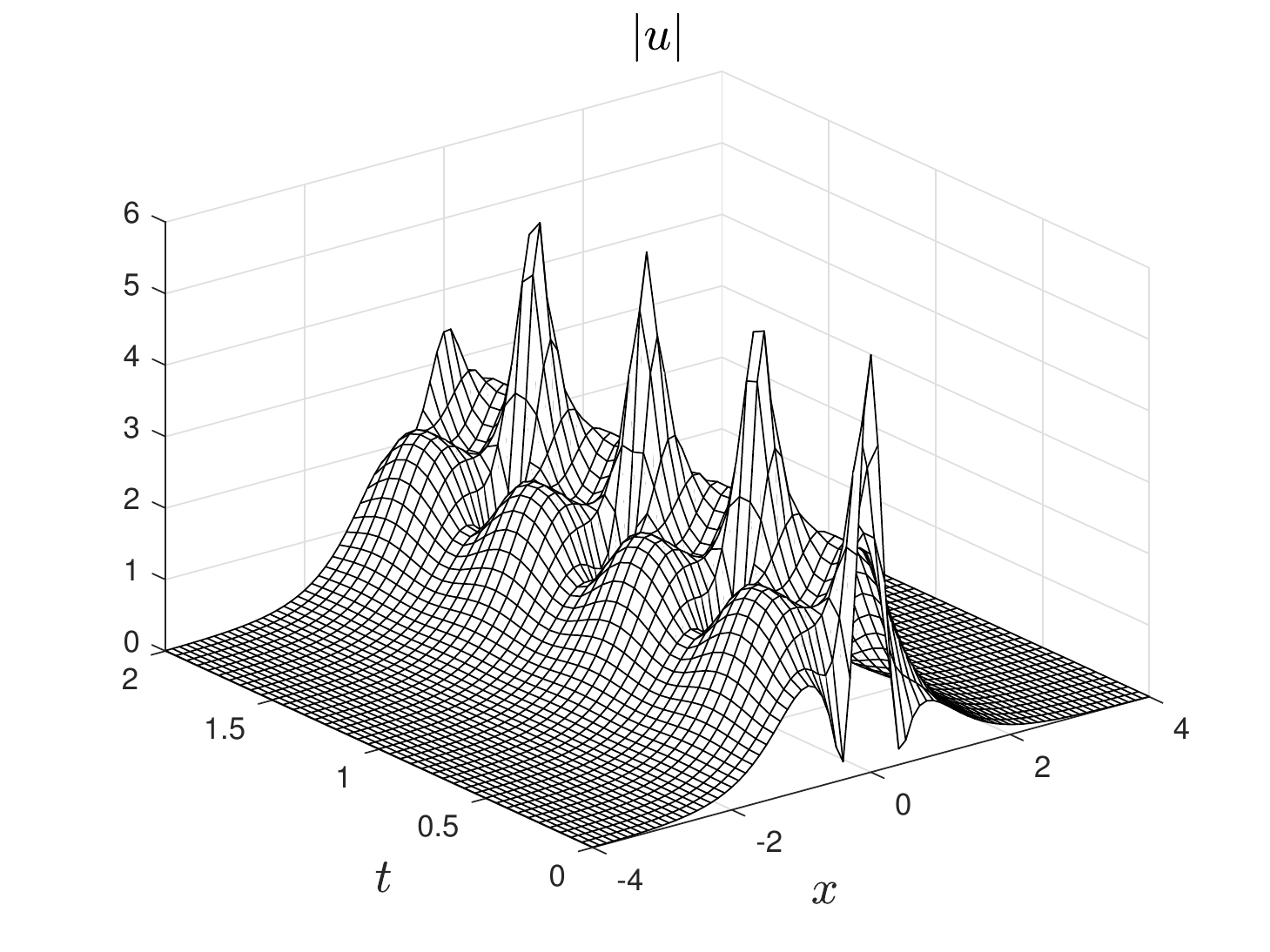}
        \caption{$|u|$}
        \label{fig:gull}
    \end{subfigure}
    ~ 
    \begin{subfigure}[b]{0.3\textwidth}
        \includegraphics[width=\textwidth]{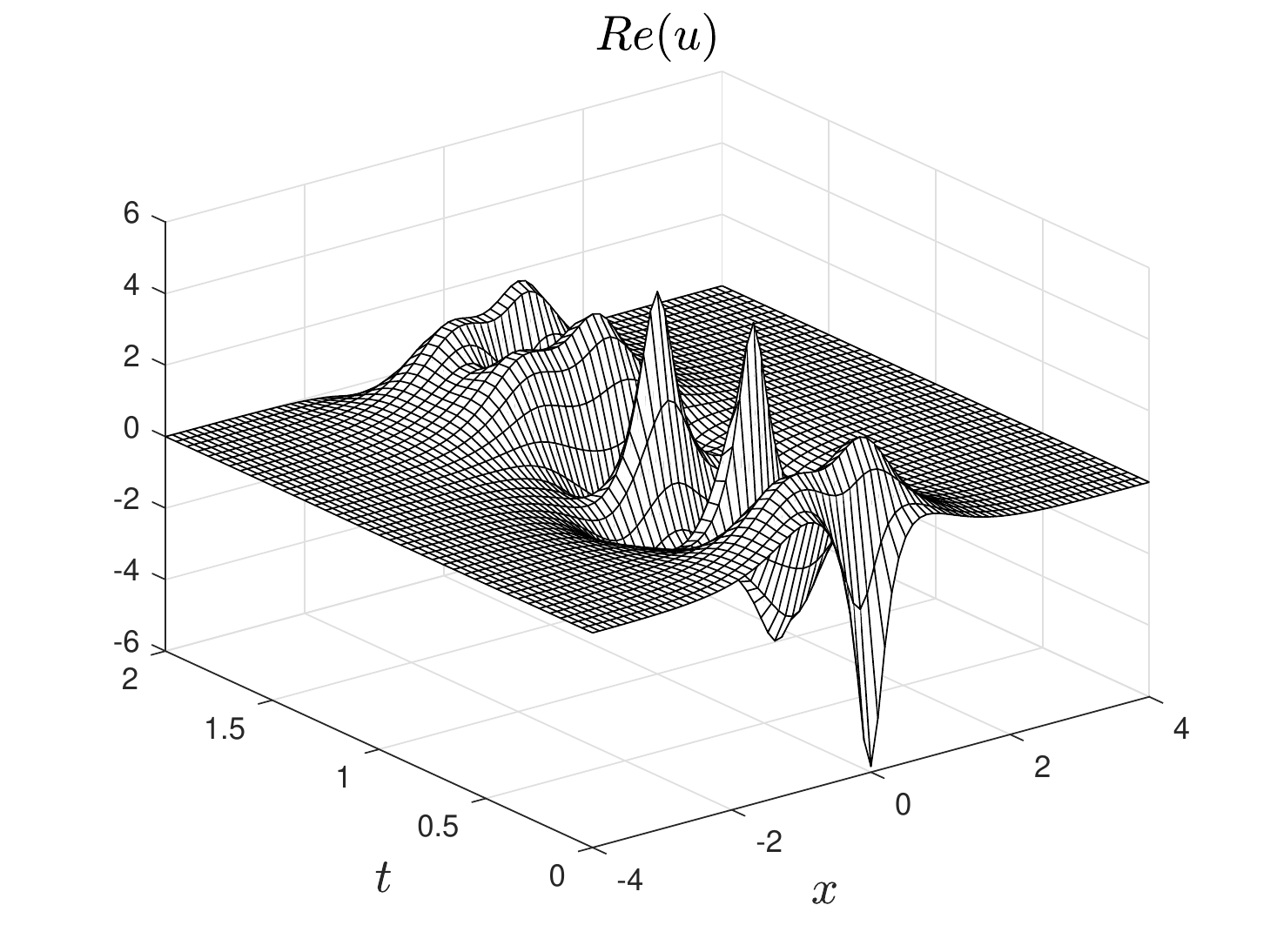}
        \caption{Re($u$)}
        \label{fig:tiger}
    \end{subfigure}
    ~ 
    \begin{subfigure}[b]{0.3\textwidth}
        \includegraphics[width=\textwidth]{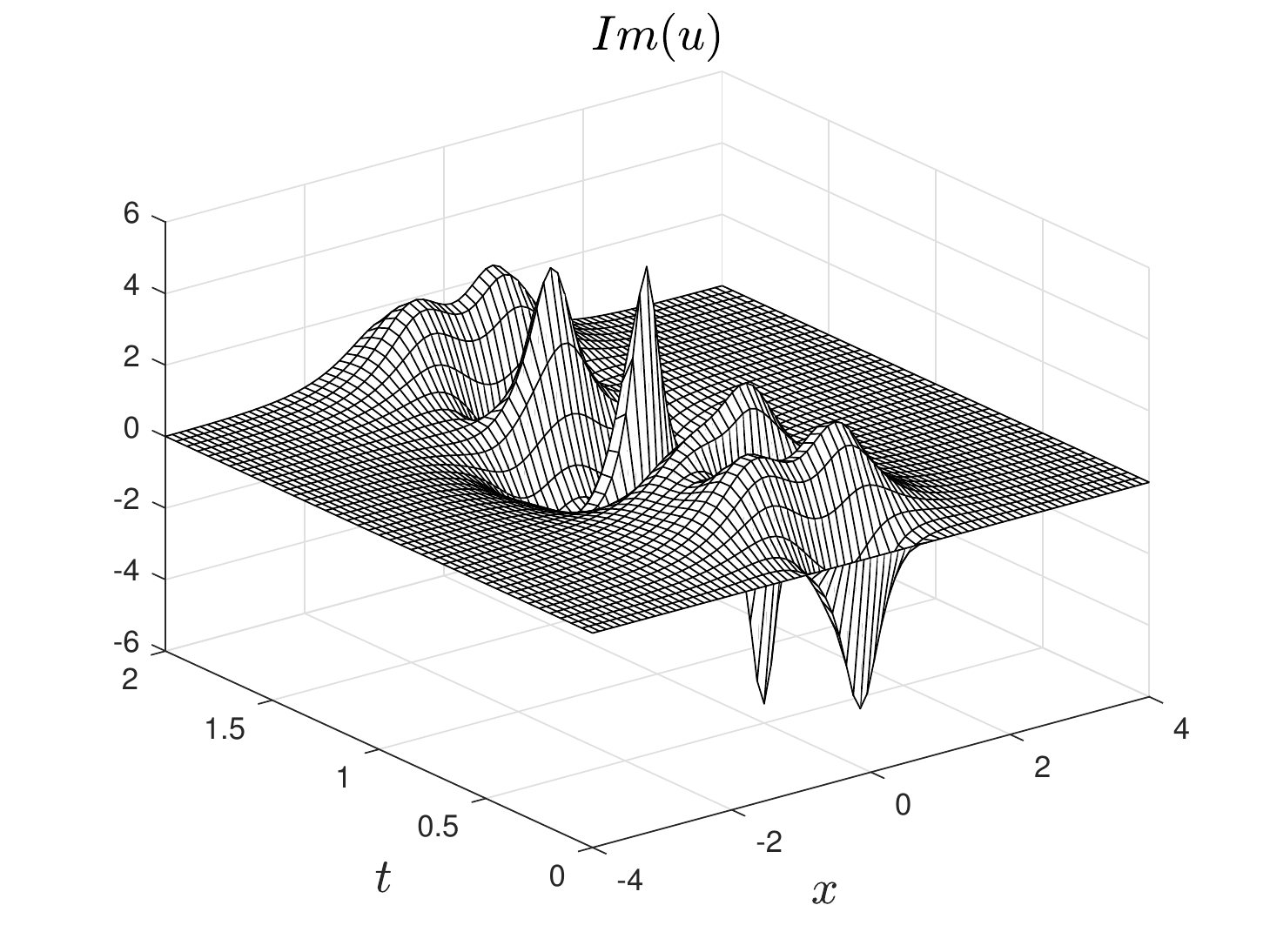}
        \caption{Im($u$)}
        \label{fig:mouse}
    \end{subfigure}
    \caption{Figures (a)-(c) show the solution to the two stationary soliton test case.  }
\label{fig:Solution}
\end{figure}

	A more difficult 
	%and novel 
	test case related to signal propagation in optical fibers and involving two interacting stationary solitons was constructed in \cite{Exact2010}. In this example we are looking for the solution $u$ of the following NLSE with a focusing nonlinearity
  	\begin{align*}
	\ci \partial_t u = - \partial_{xx} u  - 2|u|^2u \qquad \mbox{in } \mathbb{R}\times (0,T]
	\end{align*}
	and with initial value 
	\begin{align*}
	u(x,0)  =  \frac{8(9e^{-4x}+16e^{4x})-32(4e^{-2x}+9e^{2x})}{-128 + 4e^{-6x}+16e^{6x}+81e^{-2x}+64e^{2x} } & .
	\end{align*}
	As derived in \cite{Exact2010}, the exact solution is given by
\begin{align*}
u(x,t)& =\frac{8e^{4\ci t}(9e^{-4x}+16e^{4x})-32e^{16\ci t}(4e^{-2x}+9e^{2x})}{-128\cos(12t)+4e^{-6x}+16e^{6x}+81e^{-2x}+64e^{2x}}.
\end{align*}
The function $u(x,t)$ is depicted in Fig. \ref{fig:Solution} and consists of two interacting solitons.  With this it is less smooth than the single traveling soliton considered before. In fact the $L^2$-norm of its derivatives grows geometrically; already for the 9th derivative the size of the $L^2$-norm is of order $10^{11}$. However, as in the previous example, we have an exponential decay with $|u|\rightarrow 0$ for $|x| \rightarrow \infty$. Since $u(x,t)$ is stationary and periodic in time, it makes for a good long time test. For this example we also include the popular Strang splitting spectral method, henceforth abbreviated SP2 \cite{Spectral}. With this we can make a brief cross-comparison of finite element and spectral methods.

Due to the exponential decay, we restrict our computations to a finite computational domain of size $[-20,20]\times (0,T]$ and prescribe a homogeneous  Dirichlet boundary condition on both end of the spatial interval. Again, we make a numerical comparison of all five methods, where the initial value is the nodal interpolation of $u(x,0)$ in $S_h$. All the computations presented in this section are run on an Intel Core i7-6700 CPU with 3.40GHz$\times$8.

\subsubsection{Short time behaviour}
In Table \ref{rho_table} we present convergence of the density in $L^1$-norm for the mesh size $h = 40/51200$. For the same spatial discretization, $H^1$-norm convergence results are presented in Table \ref{h1_table} . From these tables we can make an important observation in terms of the convergence rate of the numerical methods. None of the methods converges right away, but the regime of an asymptotic convergence is entered at different time-steps sizes. For instance, the symplectic IM-FEM starts to show the quadratic convergence (in time) from around $\tau=2^{-8}$, the energy-conserving CN-FEM and the Two-Step FEM from around $\tau=2^{-7}$, the energy-conserving RE-FEM from around $\tau=2^{-6}$ (and even $\tau=2^{-5}$ for the $H^1$-error) and the LCN-FEM from around $\tau=2^{-8}/2^{-9}$. Due to this early start into the asymptotics, the RE-FEM scheme performs best and achieves the smallest errors. The CN-FEM is close behind, however has a higher computational complexity due to the nonlinear time-stepping.

For very large time-steps, e.g. $\tau = 1/8$, the nonlinear methods (IM-FEM and CN-FEM) become extremely ill-conditioned and all accuracy is lost. This is related to the focusing nonlinearity (i.e. $\beta<0$) and we could not observe similar effects in the defocusing case (i.e. $\beta\ge0$).

The RE-FEM performs the best and the energy conservative methods outperform by far the symplectic IM-FEM. However, taking the spectral method into account we find that it outperforms the RE-FEM. The SP2 is, for this spatial discretization ($h=40/51200$), roughly 12 times faster than the linear FE methods. The table is adjusted so that the computational times of the SP2 are comparable to those of the linear methods.

\begin{table}
	\scriptsize
	\centering
	\begin{tabular}{l |c| c | c| c| c ?  l|c|}
		
		\multicolumn{8}{c}{$||\rho-\rho_{h,\tau}||_{L^1}$}\\
		
		$\tau$	 & \multicolumn{1}{c}{IM-FEM} & \multicolumn{1}{c}{CN-FEM} & \multicolumn{1}{c}{RE-FEM} & \multicolumn{1}{c}{LCN-FEM} & \multicolumn{1}{c ?}{TwoStepFEM} & $\tau$ & \multicolumn{1}{c|}{SP2} \\
		\hline 
		$2^{-3 }$	& 2.684e14      & 3.368e12     & 12.068	  &  21.033  & 11.757  & $2^{-7}$   & 11.582 \\
		$2^{-4 }$	& 23.911       &22.546 		    & 13.774  & 18.639  &  11.0617 & $2^{-8}$  & 6.2246 \\
		$2^{-5 }$	& 23.325      &19.659   		    & 7.544   &  18.145  &  19.100 &  $2^{-9}$   & 2.1756 \\
		
		$2^{-6 }$	& 15.885      & 12.738		    & 4.648    & 18.525  & 15.223 &    $2^{-10}$  & 0.5953 \\
		
		$2^{-7 }$	& 8.498       & 5.355    	   & 1.265  & \ 4.868  & 7.733 &  $2^{-11}$& 0.1523 \\
		$2^{-8 }$	& 3.356       & 1.442           & 0.319  & \ 8.109   & \ 2.301 & $2^{-12}$& 0.0383\\
		$2^{-9 }$	& 0.959       & 0.360           & 0.080  & \ 3.537  & \ 0.577  & $2^{-13}$& 0.0096 \\
		$2^{-10 }$	& 0.252      &  0.091            &0.021  & \ 0.723  &  \ 0.145   & $2^{-14}$ & 0.0024 \\
		\hline
	\end{tabular}
	\caption{\scriptsize $L^1$-errors of the density for the two stationary soliton test case at time $T=2$.} \label{rho_table}
\end{table}

\begin{table}
	\scriptsize
	\centering
	\begin{tabular}{l |c| c | c| c| c ?  l|c|}
		
		\multicolumn{8}{c}{$||u-u_{h,\tau}||_{H^1}$}\\
		
		$\tau$	 & \multicolumn{1}{c}{IM-FEM} & \multicolumn{1}{c}{CN-FEM} & \multicolumn{1}{c}{RE-FEM} & \multicolumn{1}{c}{LCN-FEM} & \multicolumn{1}{c ?}{TwoStepFEM} & $\tau$ & \multicolumn{1}{c|}{SP2} \\
		\hline 
		$2^{-3 }$	& 3.795e10	& 4.270e9		& 11.445	&24.614		&10.629  	& $2^{-7}$  & 9.8755 \\
		$2^{-4 }$	& 10.368	& 11.108		& 15.194  	&22.414 	&11.157	& $2^{-8}$  & 4.6047\\
		$2^{-5 }$	& 9.544    & 11.474			& 13.620    &29.591  	&10.877 	& $2^{-9}$  & 1.7293 \\
		$2^{-6 }$	& 10.400    & 9.935		& 3.768		&36.725  	&10.542 	&$2^{-10}$  & 0.4877 \\
		$2^{-7 }$	& 6.175    & 4.662			& 1.022		&9.990 	&6.601		&$2^{-11}$	& 0.1258 \\
		$2^{-8 }$	& 2.704		& 1.261			& 0.260		&12.261  	&2.043 		&$2^{-12}$	& 0.0317\\
		$2^{-9 }$	& 0.806     & 0.315         & 0.066 	&3.277 	&0.514		& $2^{-13}$ & 0.0079 \\
		$2^{-10 }$	& 0.214     &  0.080        & 0.018  & 0.632  &0.129   & $2^{-14}$ & 0.0020 \\
		\hline
	\end{tabular}
	\caption{\scriptsize $H^1$-errors for the two stationary soliton test case at time $T=2$.} \label{h1_table}
\end{table}

\subsubsection{Long time behaviour}

\begin{figure}[h!] 
\begin{subfigure}[b]{0.5\textwidth}
\includegraphics[scale = 0.55]{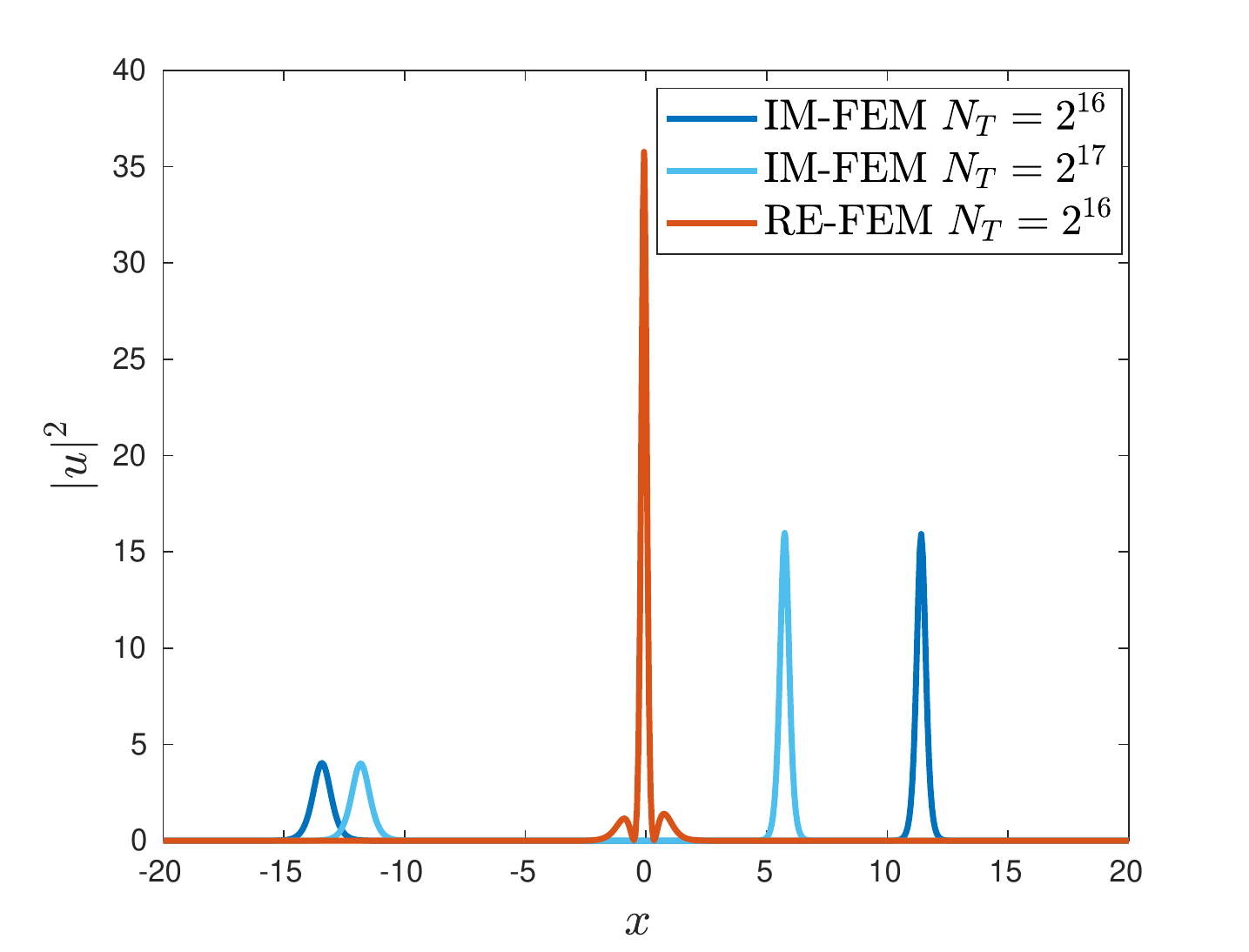}
\caption{\scriptsize Density for the RE-FEM and the IM-FEM at time $T=200$. $N_T$ denotes number of time-steps.}\label{IMvsRE}
\end{subfigure}
\begin{subfigure}[b]{0.5\textwidth}
\includegraphics[scale = 0.5]{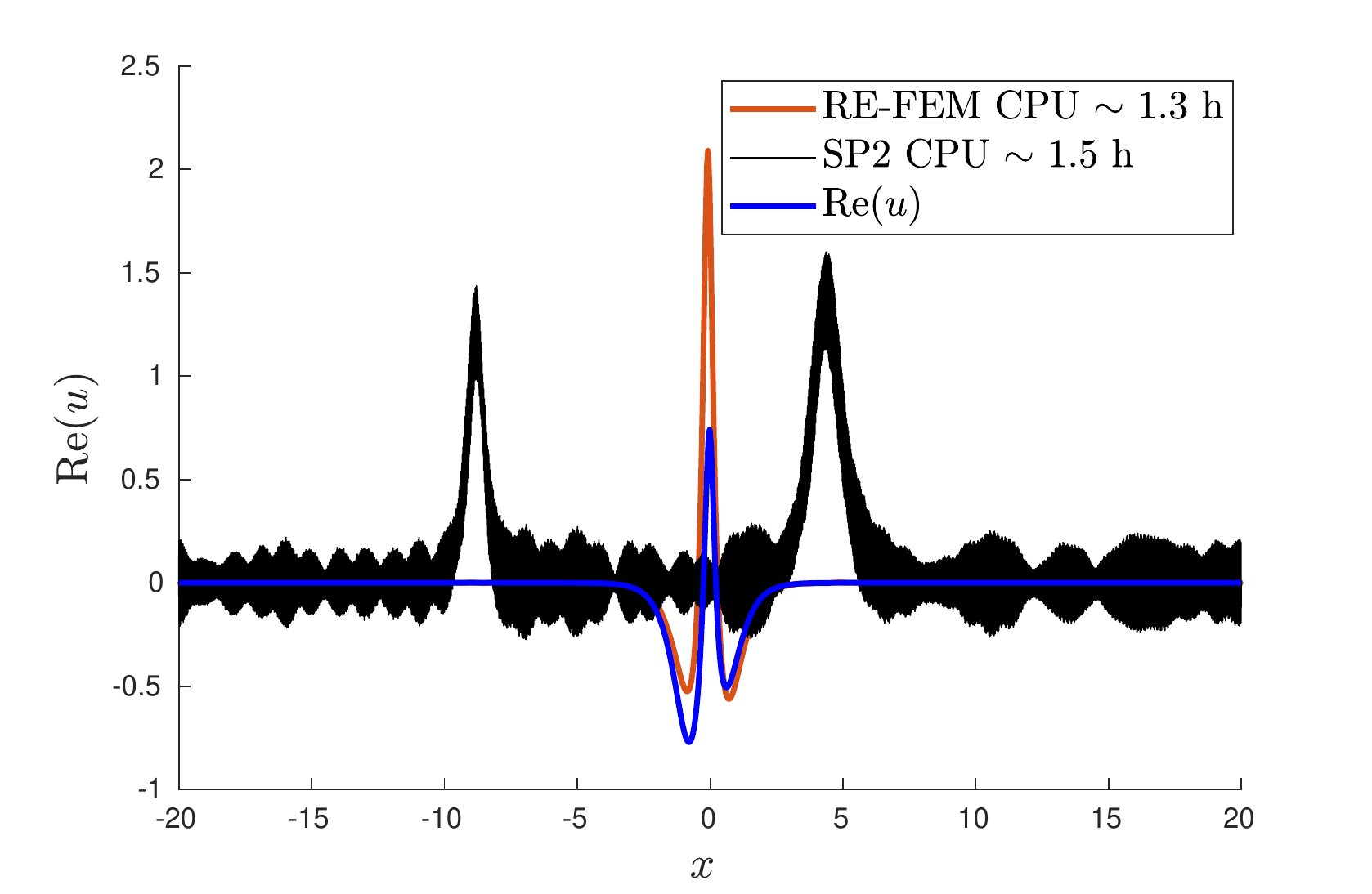}
\caption{\scriptsize Real part of solution for SP2 and RE-FEM at time $T=200$ and for comparable computational times.}\label{SP2vsRE}
\end{subfigure}
\begin{subfigure}[b]{0.5\textwidth}
\includegraphics[scale = 0.55]{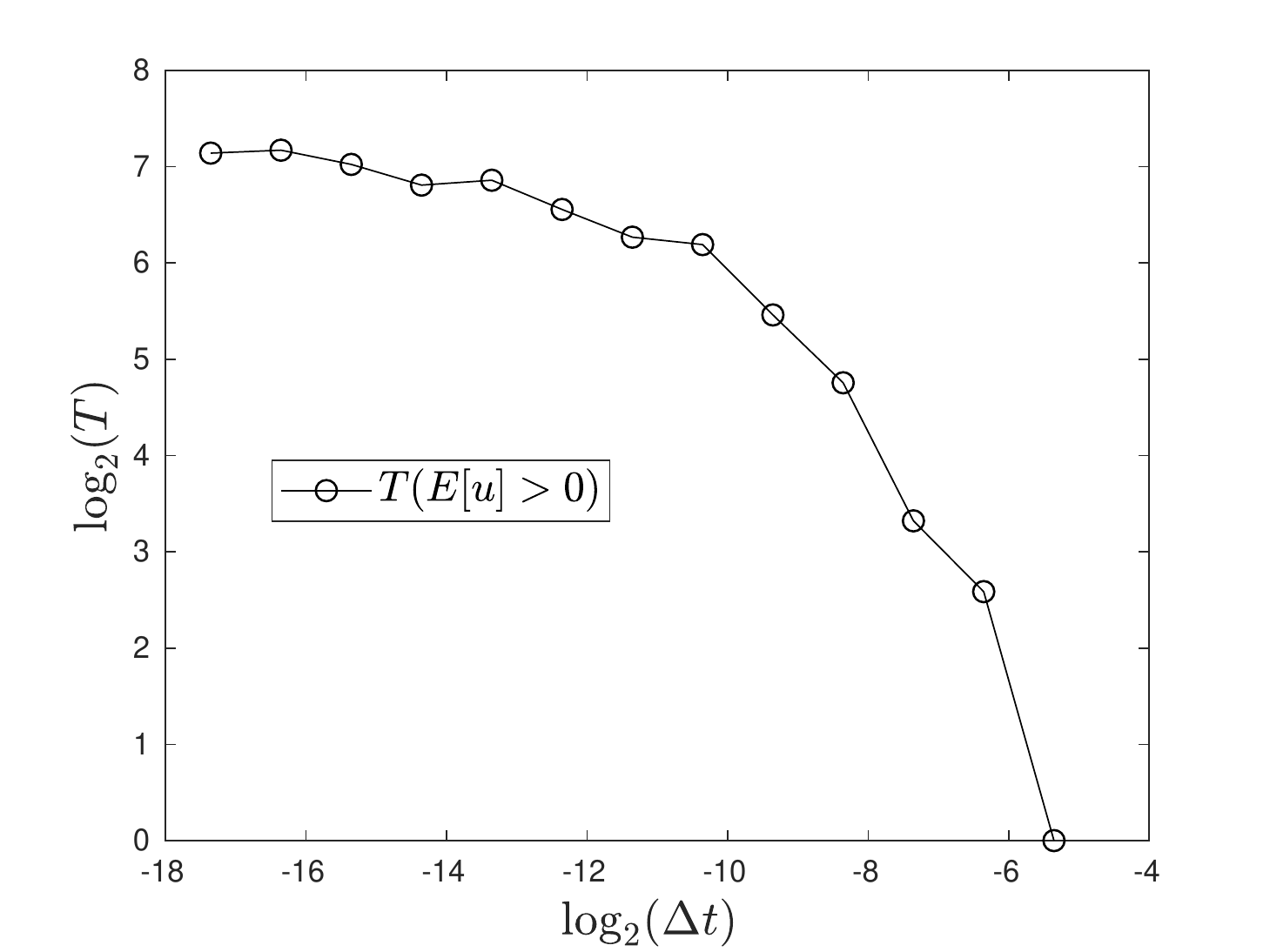}
\caption{\scriptsize Time at which the energy of the SP2 becomes positive. ($E[u] = -48$)}\label{Texp}
\end{subfigure}
\begin{subfigure}[b]{0.5\textwidth}
\includegraphics[scale = 0.5]{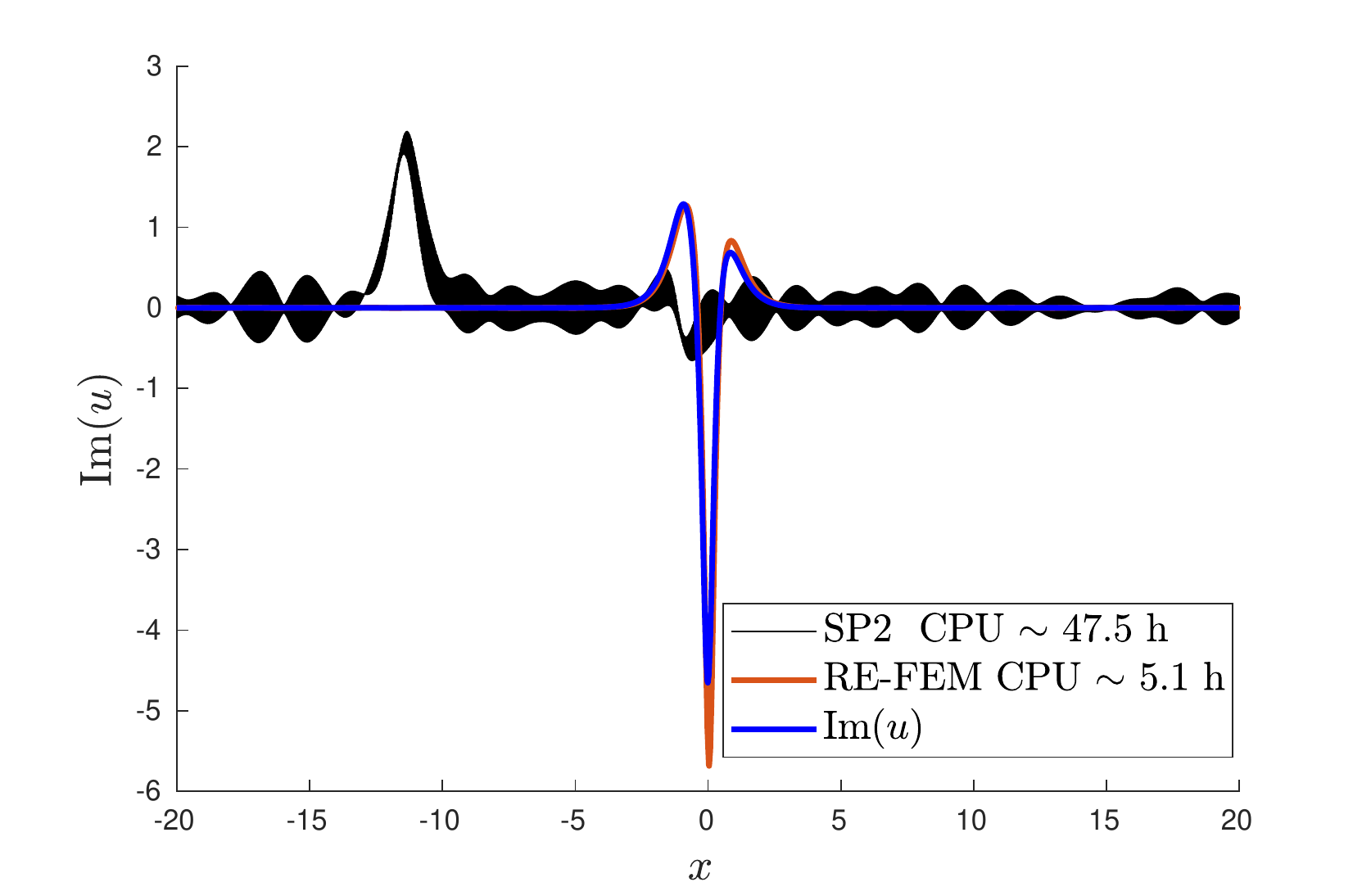}
\caption{\scriptsize Imaginary part of solution for SP2 and RE-FEM at final time $T=500$.}\label{SP2vsRE_T500}
\end{subfigure}
\caption{Long time behaviour for Model Problem \ref{StationarySoliton}.}
\end{figure}

It is often put forward that symplectic methods perform well for long time computations of Hamiltonian systems. However, for this test case we find that energy conservation is preferable. Figure \ref{IMvsRE} shows the density at $T=200$ of the IM-FEM using $2^{16}$ and $2^{17}$ time-steps, as a reference the density of the RE-FEM using $2^{16}$ time-steps is also plotted. Preserving the symplectic structure results in two separate solitons, despite the energy only varying by 0.1\% when $N_T = 2^{17}$.  In contrast, the RE-FEM captures the correct behaviour. 
	
 For these large final times, the SP2 becomes extremely prone to blow up in energy. In Figure \ref{SP2vsRE} we see the result of 1.5 hours of computation for the SP2 and the RE-FEM, corresponding to $2^{20}$ and $2^{16}$ time-steps respectively. The energy of the SP2 is almost $6\cdot 10^6$. An extreme example of this is shown in Fig. \ref{SP2vsRE_T500}, where the imaginary part is plotted at $T=500$, the SP2 uses almost 50h of CPU-time corresponding to $2^{25}$ time-steps whereas the RE-FEM uses 5h of CPU-time and $2^{18}$ time-steps. More importantly, as opposed to shorter final times, taking smaller time-steps with the SP2 only defers the blow up by an inconsiderable time. This is shown in Figure \ref{Texp}, where the time at which $E[u_{h,\tau}]>0$ for the SP2 is plotted versus time-step size ($E[u]=-48$). Each data point is a halving of the time-step size, initially halving the time-step size results in much longer energy stability. However, this diminishes quickly; even with $\Delta t = 200/2^{25}\approx 2^{-18}$, the SP2 does not compute reasonable solutions beyond $T=150$.

\subsection{Condensate in an optical lattice} 

We conclude the experiments in 1D with a test case that describes a Bose-Einstein condensate in an optical lattice. This test case demonstrates how the LCN-FEM and the Two-Step FEM can develop an instability resulting from a blow-up of the energy. As the energy remains (almost) conserved for IM-FEM, CN-FEM and RE-FEM, these methods do not suffer from such an instability. 

  For the space interval $\Omega:=[-16,16]$ we seek the solution $u$ to
\begin{align*} 
\ci \partial_t u &=-\frac{1}{2}\partial_{xx} u  +V \hspace{2pt} u+\beta \hspace{2pt} |u|^2u    
\qquad \mbox{in } \Omega,
\end{align*}
with Dirichlet boundary condition $u(-16,t)=u(16,t) = 0$ for $t\ge 0$ and initial condition $u(\cdot,0)=\uzero$ in $\Omega$. The final time is varied. The problem involves a defocusing nonlinearity with $\beta=1000$. The potential describes the combination of a harmonic confining potential together with an optical lattice and is given by
\begin{gather*} 
V(x) = (\gamma_x \hspace{2pt} x)^2+500\sin(x\pi/4)^2, \qquad \mbox{where the trapping frequency is } \gamma_x = \frac{1}{2}.
\end{gather*}
The initial state $\uzero(x)$ is selected as the ground state associated with the trapping frequency $\gamma_x = 1$. This ground state is accurately computed using the inverse iteration method for eigenvalue problems with eigenvector nonlinearities presented in \cite{Elias}. Decreasing the strength of the harmonic potential for $t>0$ (as described above) induces a small periodic motion of the condensate which, in our case, is concentrated around the lattice points $[-8, -4, 0, 4, 8]$. The density of the condensate remains almost constant. 

\begin{figure}[h!]

\centering

    \begin{subfigure}[b]{0.45\textwidth}
        \includegraphics[width=1.1\textwidth]{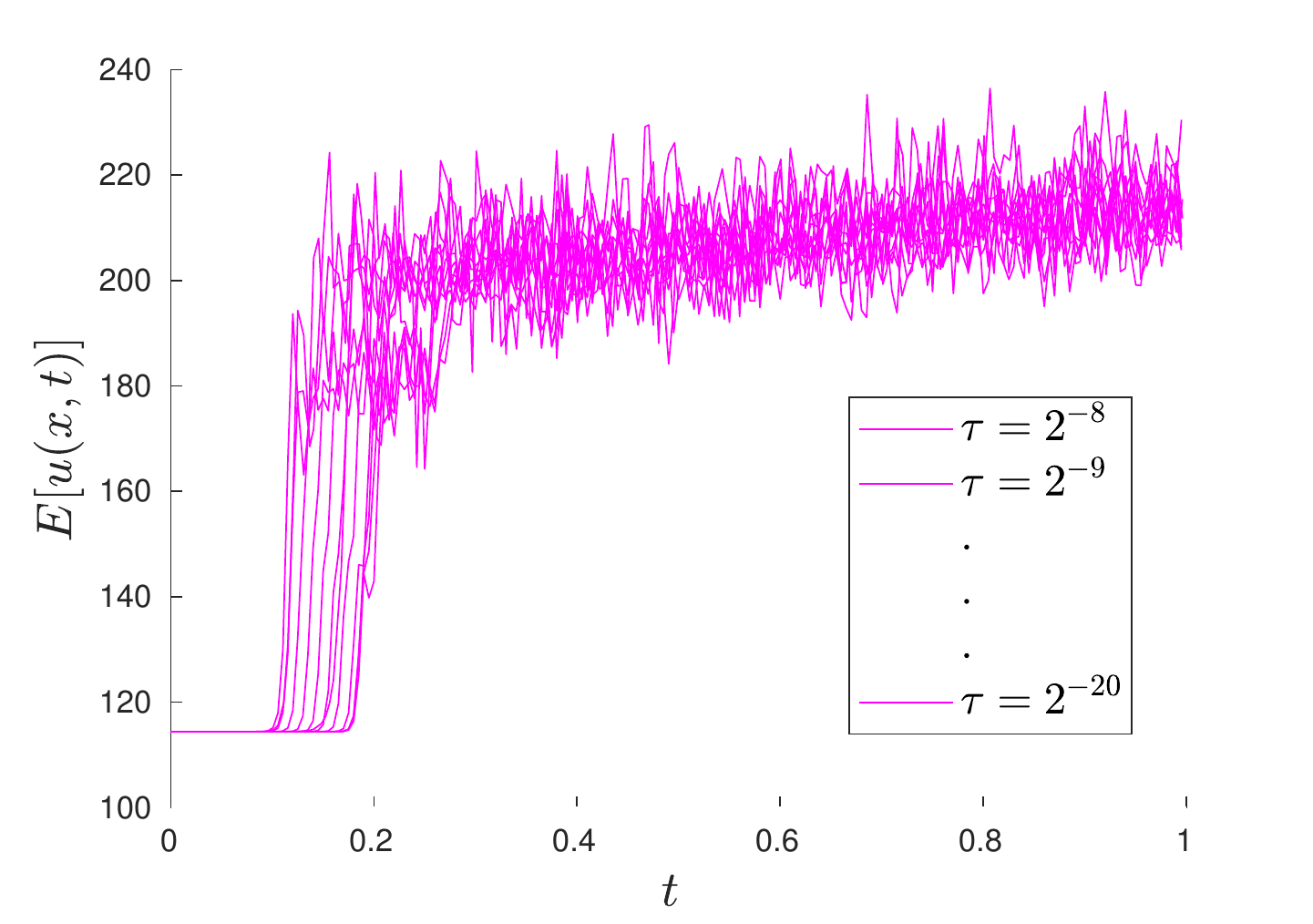} 
        \caption{\scriptsize Energy of the Two-Step FEM for the mesh size $h = 24/800$.} \label{Two_Step_Energy}
    \end{subfigure}
    ~ %add desired spacing between images, e. g. ~, \quad, \qquad, \hfill etc. 
    \begin{subfigure}[b]{0.45\textwidth}
        \includegraphics[width=1.1\textwidth]{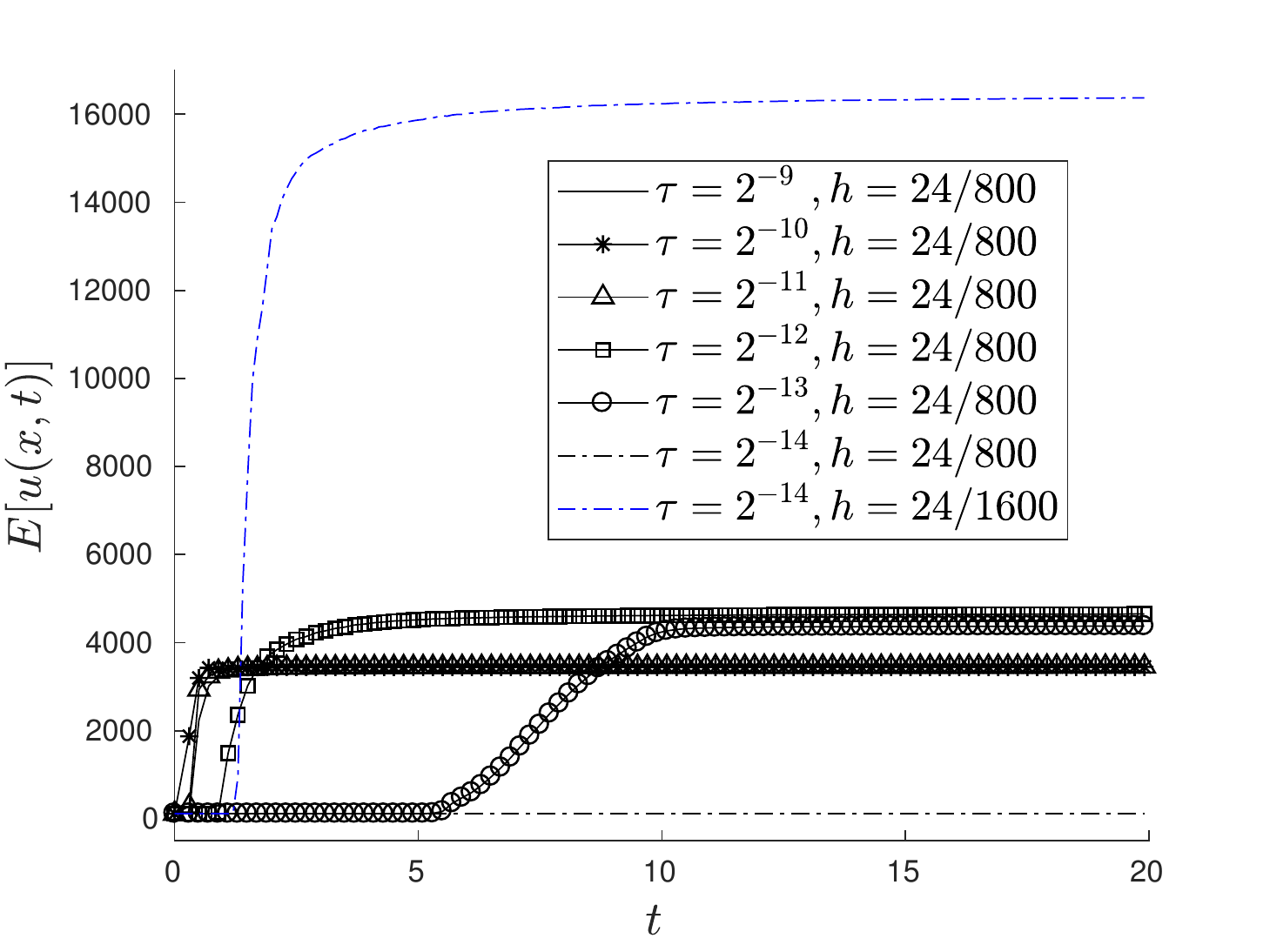} 
        \caption{\scriptsize Energy of the LCN-FEM for $h = 24/800$ in black and for $h=24/1600$ in blue. } \label{Wang_Energy}
    \end{subfigure}
    ~ %add desired spacing between images, e. g. ~, \quad, \qquad, \hfill etc. 
    \caption{Energy evolution for the Two-Step FEM and the LCN-FEM for the time discretizations $\tau = 2^{-8},2^{-9},\dots,2^{-20}$, note the different time scales. The Two-Step FEM is highly unstable for all ratios of $h$ and $\tau$. The LCN-FEM can be highly unstable depending on the ratio of $h$ and $\tau$.  }     
      \label{Instability}
\end{figure}

However, the setup is very sensitive to small perturbations of the energy and causes instabilities for the two non-symplectic / non-energy-conservative methods. This is stressed by Fig. \ref{Two_Step_Energy}, where the energy is plotted versus time for the Two-Step FEM  with spatial resolution $h = 24/800$ and time-step sizes $\tau = 2^{-8}, 2^{-9}, \dots, 2^{-20}$. For the LCN-FEM, a similar plot but for a larger final time is presented in Fig. \ref{Wang_Energy}. The spatial discretization is as before $h=24/800$ and $\tau$ ranges from $2^{-9}$ to $2^{-14}$. Energy blow-up occurs at some time for both methods of which the LCN-FEM blows up more violently. However, the most striking feature is that refining the time-step barely improves the stability of the Two-Step FEM, i.e. the blow-up is hardly delayed. On the contrary, refining the time-step greatly delays the blow-up of the LCN-FEM. Furthermore, we note that for the LCN-FEM this blow-up is coupled with the mesh size. This is also illustrated in Fig. \ref{Wang_Energy}, as changing mesh size from $h=24/800$ (black dashed line) to $h=24/1600$ (blue dashed line) while keeping $\tau = 2^{-14}$ changes the method from being stable at least till $t = 20$ to energy blow-up occurring at $t = 1.1$. Some further computations, not included for sake of brevity, show that this coupling behaves rather erratically. Finding the most stable ratio may therefore prove difficult.

Lastly we present errors and CPU-times for the five methods in Table \ref{Accuracy_CPU}. It is clear that the linear methods, with the exception of the Two-Step FEM, achieve lower errors than the nonlinear methods for the same computational times. It is also noted that the RE-FEM and CN-FEM produce almost identical solution, which is explained by the fact that the change in density is small and that for constant density the CN-FEM and the RE-FEM are equivalent. As mentioned, the Two-Step FEM does not perform well compared to the other methods due to its critical instability. Therefore and for the sake of brevity, it will henceforth not be considered in the remaining experiments.

\begin{table}[h!]
	\centering
	\small
	\begin{tabular}{ c | c  c  c  c  c | }	
		\multicolumn{6}{c}{CPU [s] times 1D}\\
		\cline{2-6}
		$N$	 & \multicolumn{1}{c}{IM-FEM} & \multicolumn{1}{c}{CN-FEM} & \multicolumn{1}{c}{RE-FEM} & \multicolumn{1}{c}{LCN-FEM} & \multicolumn{1}{c|}{TwoStep FEM} \\ \hline
		1600 & 47.92 & 46.58  & 5.78  & 5.40 & 5.50 \\
		3200 & 96.59 & 92.57 & 11.52 & 10.75 &  10.91 \\ 
		6400 & 197.89 & 192.53 & 23.38 & 21.83 &  22.18  \\
		\hline
	\end{tabular}
	\caption{Average CPU-times over five runs, the number of time-steps is 1024 and $N$ denotes the degrees of freedom. The computational complexity is the same for the linear methods and, in 1D, roughly a tenth of that of the linear methods.
	} \label{CPU1D}
\end{table}

\begin{table}[h!]
	\begin{center}
		\small

	\begin{tabular}{l | c c  |  c c |  c c|  c c|  c c| }
	\multicolumn{11}{c}{$||u-u_{h,\tau}||_{H^1}$, CPU [s]} \\ \cline{2-11}
		$\tau$	 & \multicolumn{2}{c}{IM-FEM} & \multicolumn{2}{c}{CN-FEM} & \multicolumn{2}{c}{RE-FEM} & \multicolumn{2}{c|}{LCN-FEM} & \multicolumn{2}{c|}{Two-Step FEM}  \\ 
			\cline{1-3} \cline{3-5} \cline{5-8} \cline{8-11} 
	   $2^{-9}$& 2.401  & 1220 &  3.045 & 1130 & 3.045 & 130 & 95.175 & 130 & 9.338 & 140  \\
	$2^{-10}$	& 2.474  & 1960  & 0.980 & 1800 & 0.979 & 210 & 38.324 & 210 & 9.148 & 240 \\
	$2^{-11}$	& 0.710  & 3920  & 0.249 & 3540 & 0.249 & 370 & 1.948  & 370 & 8.942 & 450 \\
	$2^{-12}$	& 0.179  & 7800  & 0.062 & 7270 & 0.062 & 700 & 0.295 & 680 & 8.899 & 870 \\
	$2^{-13}$	& 0.044  & 10750  & 0.015 & 9800 & 0.015 & 1240 & 0.080  & 1190 & 8.683 & 1610 \\ \hline
	\end{tabular}
	\caption{$H^1$-errors and CPU-times for different time-step sizes and mesh size $h=24/51200$.}\label{Accuracy_CPU}
\end{center}
\end{table}

\section{Numerical experiments in 2D}\label{EXP2D}

\subsection{Optical Lattice}
\label{optticallattice-experiment-section}

In the next experiment we will consider an optical lattice in 2D. This example serves to show that despite being both mass-conservative and symplectic, the IM-FEM has an instability that manifests itself in the density. This instability is even more pronounced in the LCN-FEM, but fully absent in the energy-conservative methods. In the corresponding test problem, we seek $u(x,t)$ with
	\begin{gather}  \label{OpticalLattice2D}
	\begin{cases}
	\ci \partial_t u  &= -\frac{1}{2} \Delta u + V u  + 2300|u|^2u   \qquad \mbox{in } \Omega \times (0,T], \\
	u(\cdot,t)  &=  0     \qquad\hspace{112pt} \mbox{on } \partial \Omega \times (0,T],  \\
	u(\cdot,0)  &= \uzero  \qquad\hspace{107pt} \mbox{in }  \Omega.
	\end{cases}
	\end{gather}
	Here, $\Omega = [-6,6]^2$ is the computational domain and for the maximum time we selected $T = 1$.
	The potential $V  = V_{\text{o}} +V_{\text{ha}} + V_{\text{c}}$, consists of three parts: an optical lattice $V_{\text{o}}$, a harmonic confining potential $V_{\text{ha}}$ and a discontinuous confining potential frame $V_{\text{c}}$. We have
\begin{align*}
V_{\text{o}}(x_1,x_2) &= \ 787 \sum_{i=1}^2 \sin(\pi x_i/2)^2, \qquad
V_{\text{ha}}(x_1,x_2)=  \frac{1}{2} \sum_{i=1}^2 (\gamma_{x_i} x_i)^2, 
\quad \mbox{and}\\ 
V_{\text{c}}(x_1,x_2)
&= 1000(\ (|x_1|-4.5)^5\chi_{|x_1|\geq 4.5 } + (|x_2|-4.5)^5\chi_{|x_2|\geq 4.5 }.
\end{align*}
For the dynamics in the time-dependent problem we set the trapping frequencies to $\gamma_x = 4$ and $\gamma_y = 8$. The initial value $\uzero$ is the ground state with $\gamma_x=\gamma_y = 1$, i.e. it solves the eigenvalue problem
\[\lambda_0 \uzero = -\frac{1}{2}\Delta \uzero +V\uzero +2300|\uzero|^2\uzero, \]
with ground state eigenvalue (chemical potential) $\lambda_0$.
The initial state is accurately computed using the inverse iteration method for eigenvalue problems with eigenvector nonlinearities presented in \cite{Elias}. The solution to \eqref{OpticalLattice2D} approximately consists of 25 localized oscillators.

\begin{figure}[h!]
    \centering
    \begin{subfigure}[b]{0.45\textwidth}
        \includegraphics[width=1.1\textwidth]{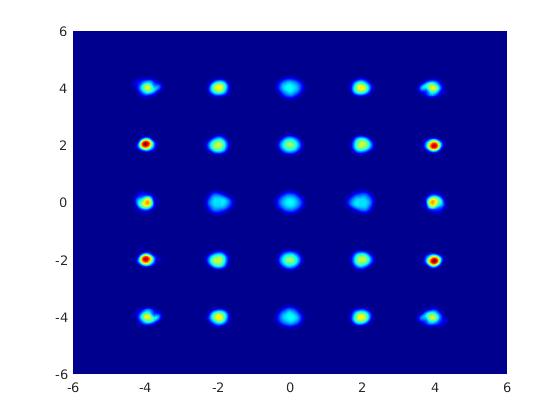}
        \caption{IM-FEM }
        \label{IM_SOL_2D-Lattice}
    \end{subfigure}
    ~ 
    \begin{subfigure}[b]{0.45\textwidth}
        \includegraphics[width=1.1\textwidth]{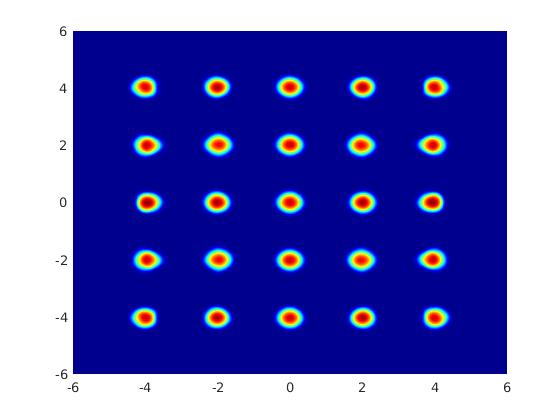}
        \caption{CN-FEM }
        \label{CN-FEM_SOL_2D-Lattice}
    \end{subfigure}
    ~
    ~

    \begin{subfigure}[b]{0.45\textwidth}
        \includegraphics[width=1.1\textwidth]{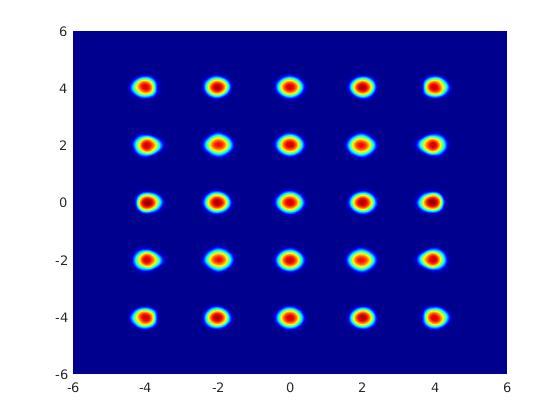}
        \caption{RE-FEM}
        \label{RE-FEM_2D-Lattice}
    \end{subfigure}
    ~ 
    \begin{subfigure}[b]{0.45\textwidth}
        \includegraphics[width=1.1\textwidth]{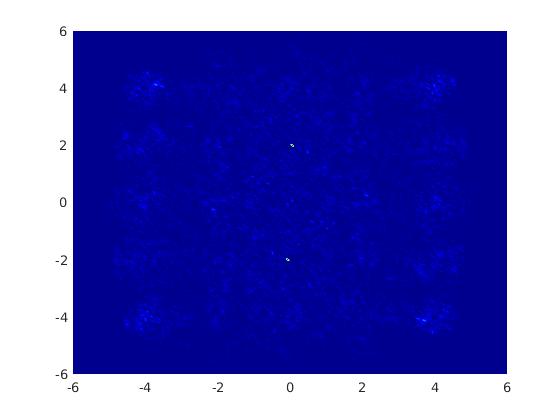}
        \caption{LCN-FEM }
        \label{LCN-FEM_2D-Lattice}
    \end{subfigure}
    ~
    \caption{Solution plots of density $|u_{h,\tau}|^2$ in Model Problem \ref{optticallattice-experiment-section} for $\tau = 2^{-8}$ and $h=0.06$.}
\end{figure}

\begin{figure}[h!]
    \centering
    \begin{subfigure}[b]{0.45\textwidth}
        \includegraphics[width=1.1\textwidth]{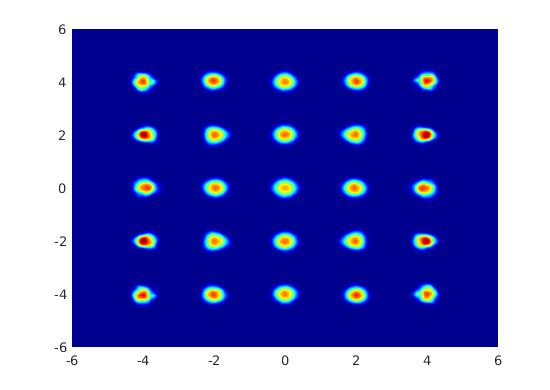}
        \caption{IM-FEM }
        \label{IM_SOL_2D-Lattice2}
    \end{subfigure}
    ~ 
    \begin{subfigure}[b]{0.45\textwidth}
        \includegraphics[width=1.1\textwidth]{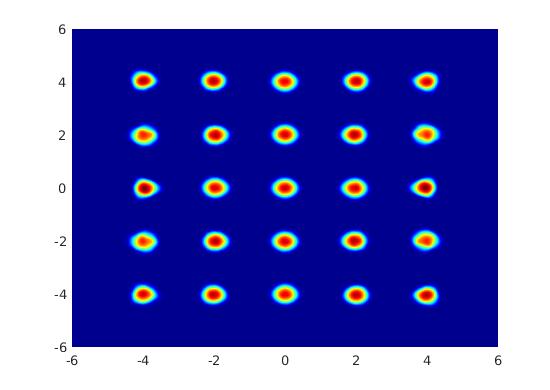}
        \caption{CN-FEM }
        \label{CN_SOL_2D-Lattice2}
    \end{subfigure}
    ~ 
    
~
    \begin{subfigure}[b]{0.45\textwidth}
        \includegraphics[width=1.1\textwidth]{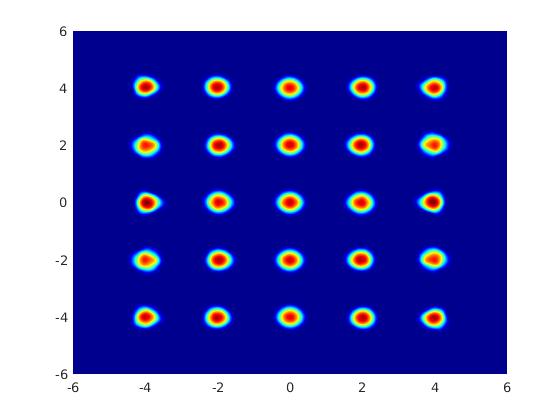}
        \caption{RE-FEM}
        \label{RE-FEM_SOL_2D-Lattice2}
    \end{subfigure}
    ~ 
    \begin{subfigure}[b]{0.45\textwidth}
        \includegraphics[width=1.1\textwidth]{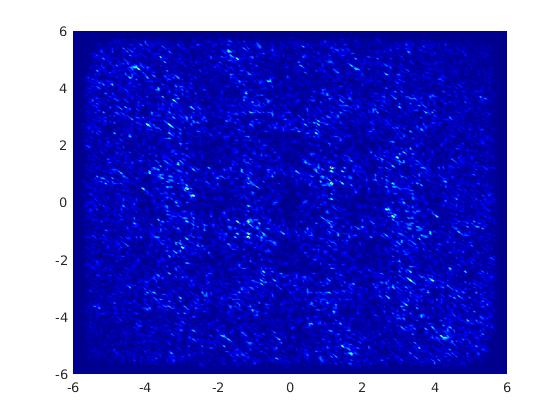}
        \caption{LCN-FEM }
        \label{LCN-FEM_SOL_2D-Lattice2}
    \end{subfigure}
    ~ 
    \caption{Solution plots of density $|u_{h,\tau}|^2$ in Model Problem \ref{optticallattice-experiment-section} for $\tau = 2^{-9}$ and $h=0.06$.}
\end{figure}

\begin{figure}[h!]
    \centering
    \begin{subfigure}[b]{0.45\textwidth}
        \includegraphics[width=1.1\textwidth]{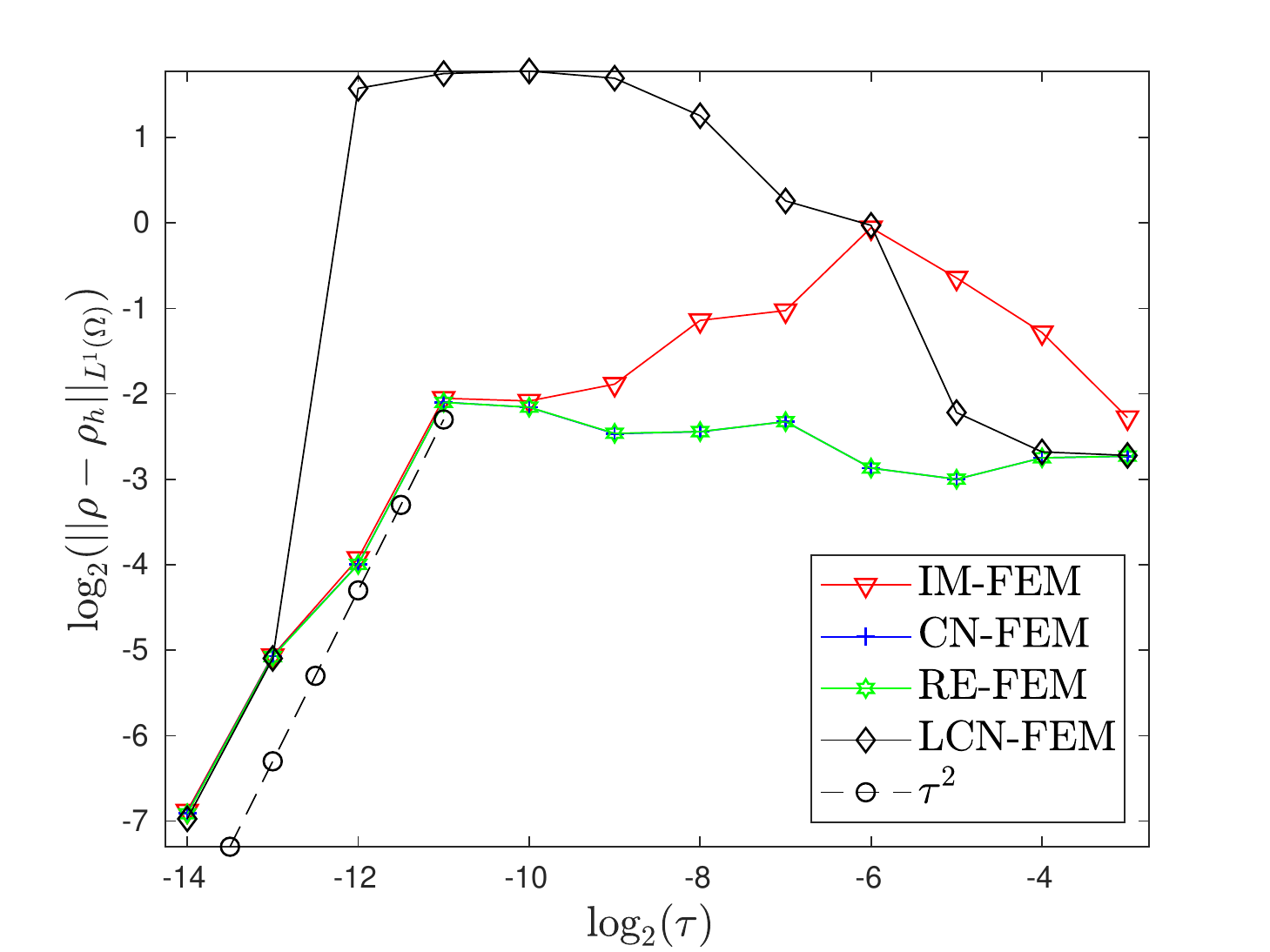}
        \caption{$||\rho-\rho_{h,\tau}||_{L^1(\Omega)}$}
        \label{SOLs_Rho}
    \end{subfigure}
    ~
    \begin{subfigure}[b]{0.45\textwidth}
        \includegraphics[width=1.1\textwidth]{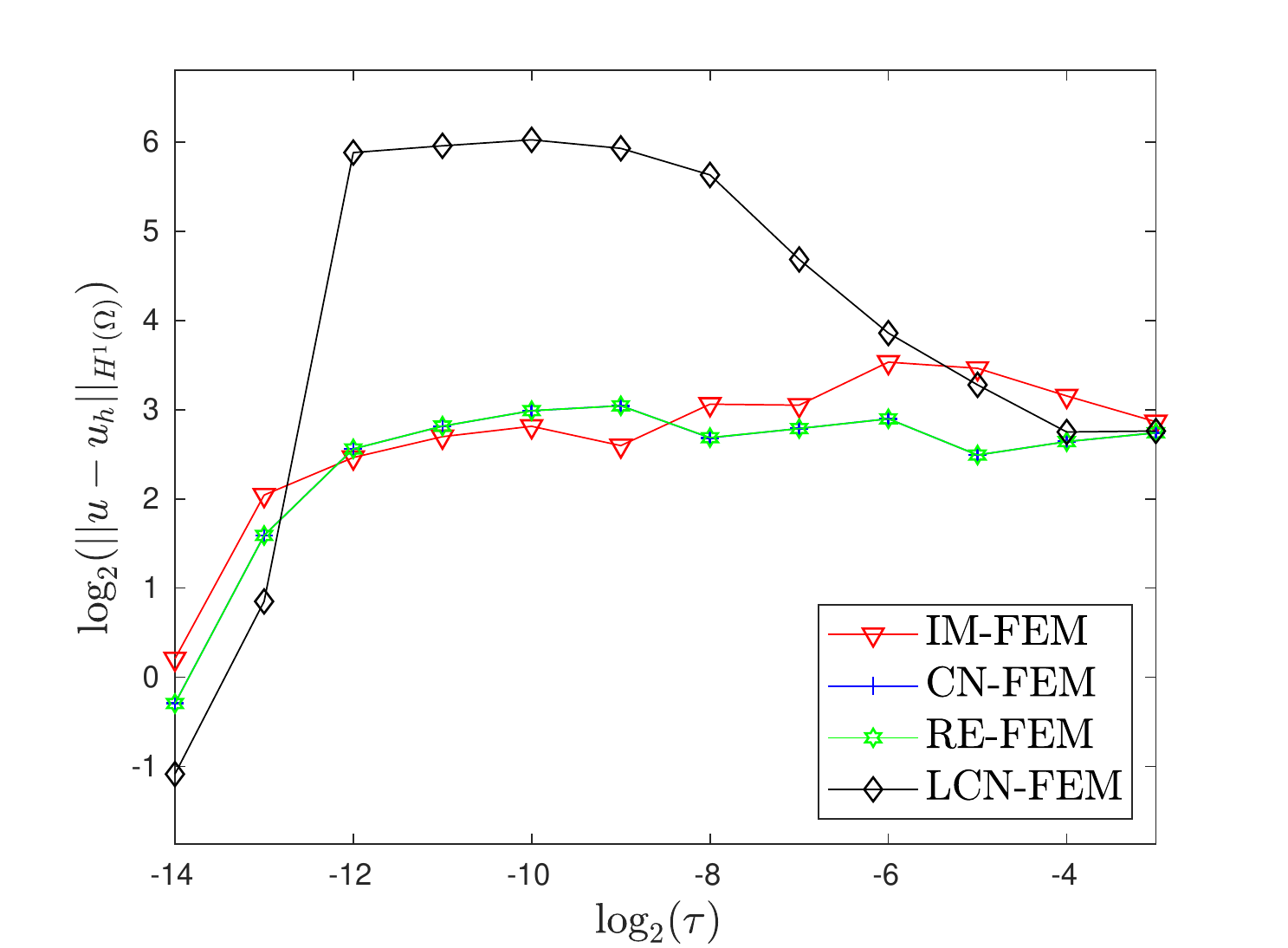}
        \caption{$||u-u_{h,\tau}||_{H^1(\Omega)}$}
        \label{SOLs_H1}
    \end{subfigure}
    ~ 
    \caption{For Model Problem \ref{optticallattice-experiment-section}, the $L^1$-errors in density and the $H^1$-errors in $u$ for the four methods versus time-step size for the spatial discretization $h=0.06$.}
\end{figure}

Again the mass-conservative LCN-FEM is prone to blow up in energy as indicated in Fig. \ref{SOLs_H1}. This results in meaningless solution plots as can be seen in Fig. \ref{LCN-FEM_2D-Lattice} and Fig. \ref{LCN-FEM_SOL_2D-Lattice2}. Fig. \ref{SOLs_Rho} shows that there is a pre-asymptotic regime and furthermore that for $h=0.06$, the LCN-FEM requires 4 times higher resolution as the other methods to enter this asymptotic regime, namely $\tau = 2^{-13}$. This is clearly coupled with $h$ as we find for $h=0.03$ that the energy blows up even for $\tau = 2^{-14}$.

\begin{figure}[h!]
 \centering
 \includegraphics[scale=0.52]{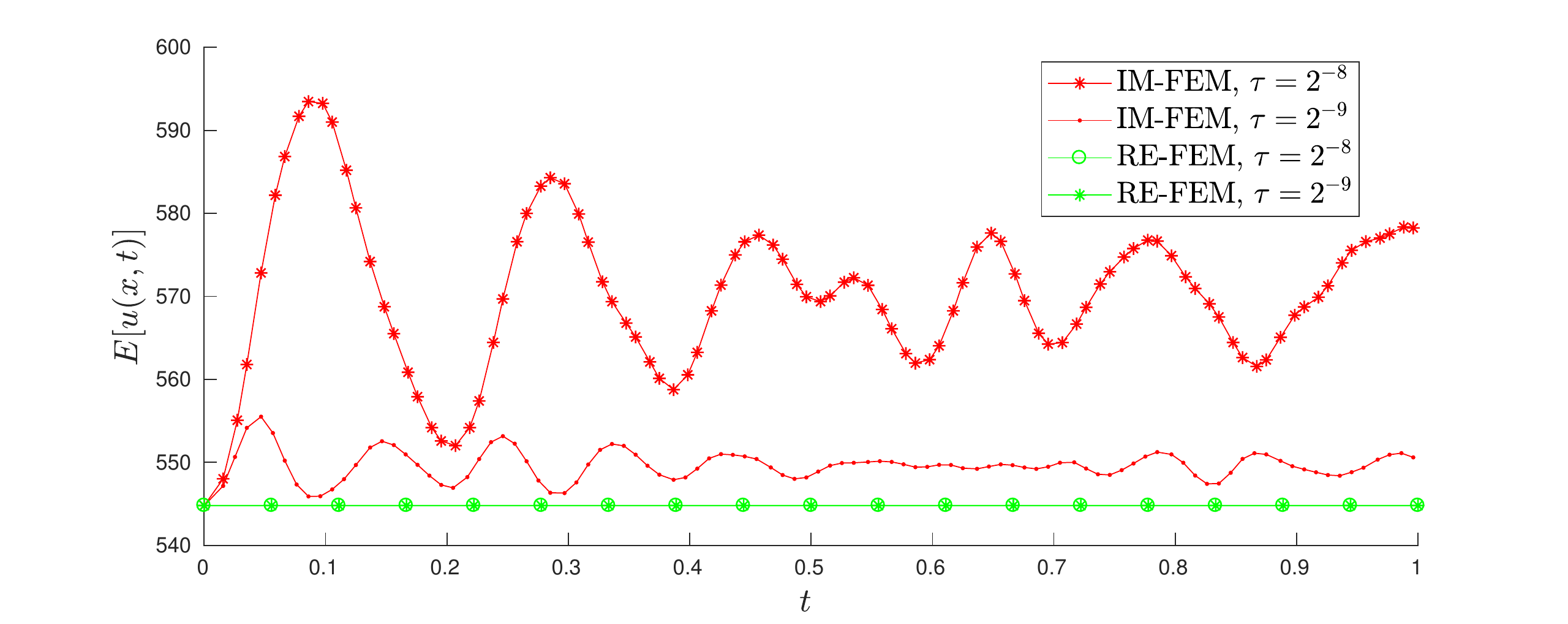}
\caption{Energy evolution of the IM-FEM and RE-FEM for two time-step sizes and the spatial discretization $h=0.06$. The corresponding solutions to Model Problem \ref{optticallattice-experiment-section} are depicted in Fig. \ref{IM_SOL_2D-Lattice}-\ref{LCN-FEM_SOL_2D-Lattice2}} \label{Energy_IM_vs_RE}
  \end{figure}

Perhaps more notably, the symplectic IM-FEM produces a deteriorated solution for coarser time-steps as the error in density blows up, cf. Fig. \ref{IM_SOL_2D-Lattice} and Fig. \ref{IM_SOL_2D-Lattice2}. Since the symplectic method does not allow a blow-up of the energy, the deterioration seems to be already triggered by small fluctuations in the energy. This guess is supported by the observation that the energy-conservative methods perform equally well and never produce completely erroneous density plots. The energy evolution of the symplectic IM-FEM is shown in Fig. \ref{Energy_IM_vs_RE}, where the energy oscillations become clearly visible.

\subsection{Rotating Bose Einstein Condensate}\label{rotatingcondensate-section}
The next experiment shows how the low errors of the LCN-FEM can be misleading as the method may still produce polluted plots.

In this experiment we consider a rotating Bose-Einstein condensate (BEC) which is typically used to investigate superfluidity on an observable scale (cf. \cite{Aft06,FSF01}). Here the superfluid character of the condensate is expressed through density singularities, so called-vortices, and it is of interest to study the dynamical behavior of such vortex patterns. In an appropriate mathematical model, we first compute the initial value as the ground state of a BEC under angular momentum rotation, i.e. we seek an $L^2$-normalized eigenfunction $\uzero \in H^1_0(\Omega)$ and corresponding minimal eigenvalue $\lambda >0$ such that 
\begin{align*}
-\frac{1}{2} \triangle \uzero + V_0 \hspace{2pt} \uzero - \boldsymbol \omega \mathcal{L} \uzero + \beta |\uzero|^2 \uzero = \lambda \uzero \qquad \mbox{in} \enspace \Omega.
\end{align*}
Here, 
$$
\mathcal{L} := - \ci \left( x \partial_y - y \partial_x \right)
$$
denotes the axis-oriented angular momentum operator and $\boldsymbol \omega$ is the angular velocity. In our experiment we selected $\Omega = (-6,6)^2$, $V_0(x,y)=\frac{1}{2}x^2+\frac{1}{2}y^2$, $\beta=100$ and $ \boldsymbol \omega = 0.8$. The eigenvalue problem was solved numerically using the Discrete Normalized Gradient Flow method proposed in \cite{BaD04}. The computational mesh is the same as the mesh used for solving the time-dependent problem afterwards. Here we simulate the situation that the stirring potential (represented by $\boldsymbol \omega \mathcal{L}$) is switched off after the condensate reaches the ground state. Due to the initial rotation the condensate still has an angular momentum and we can study the dynamical change of the vortex pattern. Modifying the trapping potential so that it models a slightly anisotropic trap, we consider the problem to find $u$ such that 
\begin{align*}
\ci \partial_t u &= -\frac{1}{2} \triangle u + V  \hspace{2pt} u + \beta |u|^2 u \qquad \mbox{in} \enspace \Omega,\\
\nonumber u&=0 \qquad \hspace{102pt} \mbox{on} \enspace \Omega
\end{align*}
with anisotropic potential $V(x,y)= 0.45 x^2 + 0.55 y^2$ and the initial condition $u(x,0)=\uzero$, where $\uzero$ is the ground state from above. The computations are carried out for maximum time $T=10$. As a reference solution in all computations we use an approximation obtained with the RE-FEM with a time-step size $\tau=2\cdot 10^{-4}$ and mesh size $h=12\cdot2^{-6}$. The mesh size will be constant for all computations so that we can focus on the error caused by the time discretization alone.

\begin{figure}[h]
 \centering
\includegraphics[width=1.0\textwidth]{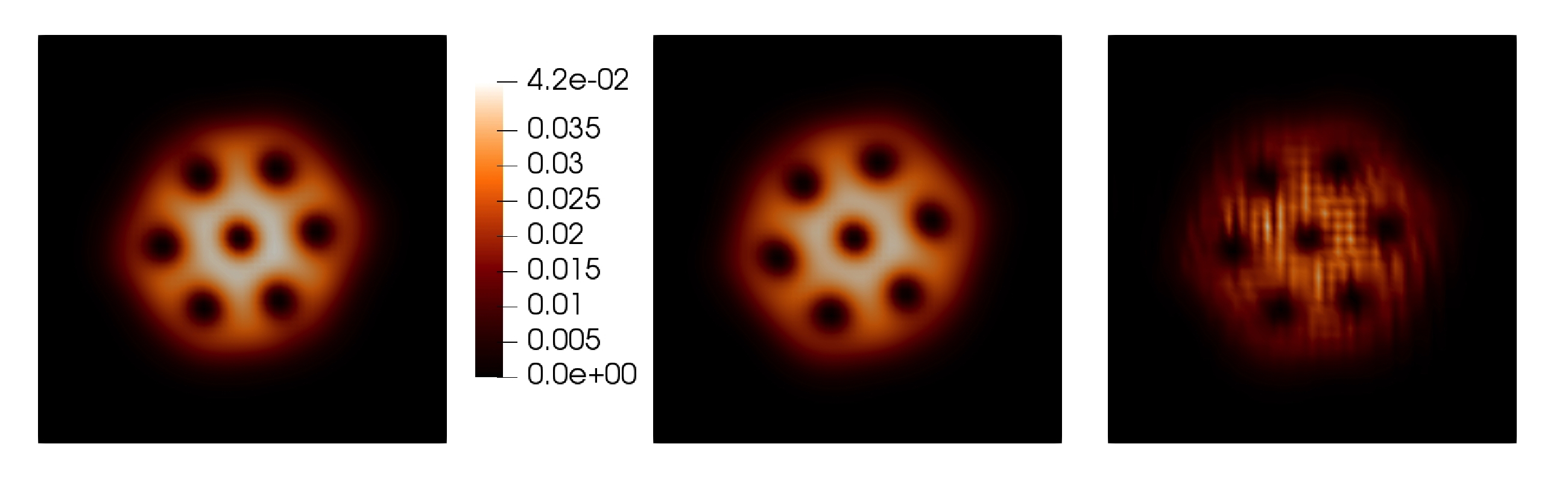}
\caption{Numerical approximations at $T=10$ and mesh resolution $h=12\cdot2^{-6}$. The left picture depicts the reference solution. The middle figure shows the solution obtained with the Crank-Nicolson FE-method for $\tau=T/256$. We note that the approximations obtained with IM-FEM and RE-FEM for the same mesh size can visually not be distinguished from the CN-FEM result. However, the approximation obtained with LCN-FEM for $\tau=T/256$ as depicted in the right figure clearly suffers from an energy pollution and looks distorted compared to the other solutions.}\label{Rotating_Density_Plots}
\end{figure}

This experiment further corroborates the general conclusions that in our experiments the energy conservative methods are slightly more accurate than the symplectic method and that the linearized Crank-Nicolson method is prone to cause unphysical effects. In particular, the LCN-FEM can produce erroneous density plots while still maintaining low errors. A striking example of this is Fig. \ref{Rotating_Density_Plots} together with the Table \ref{Rotating_Condensate_Table}. The figure shows the density plots of the LCN-FEM, CN-FEM and the RE-FEM for a fine spatial discretization and $\tau = T/256$. The LCN-FEM solution is polluted but yet the table shows that it has the lowest error in both $L^2$- and $H^1$-norm as well as the lowest error in terms of $L^1$-norm of the density. The reason for this observation is that LCN-FEM captures the rotational phase significantly better than the other methods. As for a comparison between IM-FEM, CN-FEM and RE-FEM we observe that all methods have a similar accuracy, where the errors become smallest for the energy-conserving CN-FEM. However, noting the considerably lower computational complexity of the energy-conserving RE-FEM, we can clearly identify it as the best-performing method for the test case.

\begin{figure}[h]
 \centering
\includegraphics[width=0.7\textwidth]{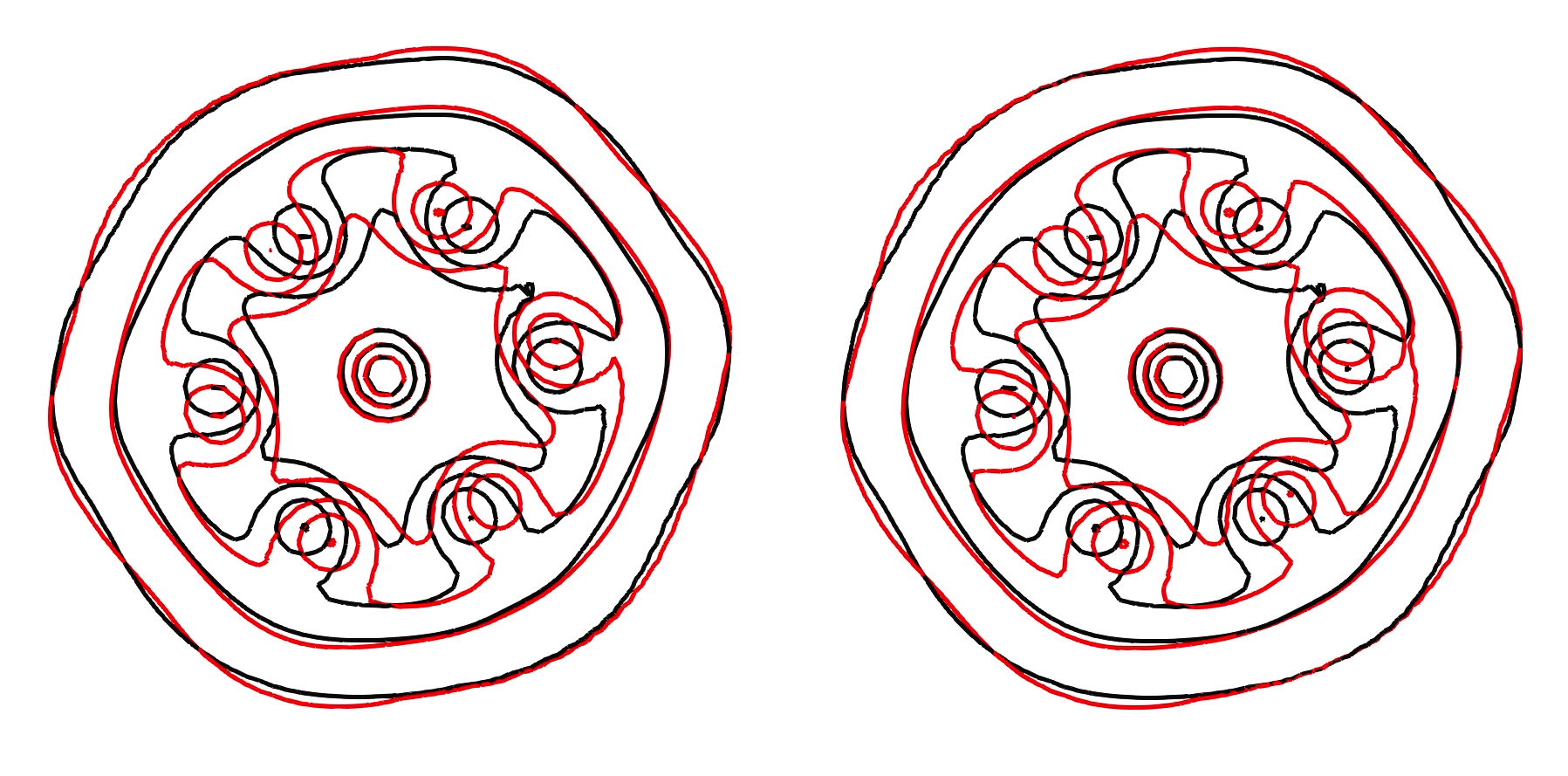}
\caption{Both figures depict a comparison between the isolines of the reference solution (black) and the isolines of a RE-FEM approximation (left, red lines) and a IM-FEM approximation (right, red), respectively. The snapshots are at maximum time $T=10$ and the approximations are computed for $\tau=T/256$ and $h=12\cdot2^{-6}$.}
\label{Isolines}
\end{figure}

It is also worth to mention that all methods suffer from a rotational phase shift compared to the exact solution. This is shown exemplarily for the IM-FEM and RE-FEM in Fig. \ref{Isolines}. The isolines between the two numerical approximations are hardly distinguishable, where the RE-FEM is only slightly better at respecting the isolines of the exact solution than the IM-FEM. One such slight difference can be seen in Fig. \ref{Isolines} for the rightmost vortex, where the IM-FEM produces a closed loop whereas the RE-FEM does not.

\begin{table}[h!]
	\begin{center}
		\tiny

	\begin{tabular}{l |l c r | l c r | l c  r| l c r| }
	\multicolumn{13}{c}{$||u-u_{h,\tau}||_{L^2}$, $||\rho-\rho_{h,\tau}||_{L^1}$, $||u-u_{h,\tau}||_{H^1}$} \\ \cline{2-13}
		$T/\tau$	 & \multicolumn{3}{c}{IM-FEM} & \multicolumn{3}{c}{CN-FEM} & \multicolumn{3}{c}{RE-FEM} & \multicolumn{3}{c|}{LCN-FEM}  \\ 
			\cline{1-3} \cline{3-5} \cline{5-8} \cline{8-13} 
		64  & 0.6590 & 0.3560 & 1.9949  &    1.7527 & 0.3594 & 3.9750 & 1.7272  & 0.3553 & 3.9466  & 1.4095  &  1.4014  & 6.6062   \\
		128 & 1.1946 & 0.3496 & 2.1666 &   1.0690  & 0.3579 & 2.8822 & 1.0779  & 0.3579 & 2.9085 & 1.4500   & 1.2261 & 8.1410   \\
		256 & 1.1855 & 0.1893 & 3.1790 &    0.8746  & 0.1714 & 2.4329 & 0.8717  & 0.1707 & 2.4257 & 0.5402   & 0.1431 & 1.6532   \\
		512 & 0.3416 & 0.0515 & 0.9392 &  0.2362  &  0.0453 &  0.6699 & 0.2373 &  0.0456 & 0.6734  & 0.6026 & 0.4097 & 7.1779 \\
		1024 & 0.0847 & 0.0126 & 0.2333 & 0.0575 & 0.0110 & 0.1636 & 0.0599 & 0.0116 &  0.1706 & 0.0400 & 0.0118 & 0.2426 \\
		2048 & 0.0192 & 0.0027 & 0.0523 & 0.0110 & 0.0023 & 0.0348 & 0.0150 & 0.0029 & 0.0427 & 0.0111 & 0.0013 & 0.0239 \\
		\hline
	\end{tabular}
	\caption{\scriptsize $L^2$ and $H^1$-errors and $L^1$-errors of the density for the rotating condensate experiment (i.e. Model Problem \ref{rotatingcondensate-section}) at time $T=10$ and for the mesh size $h=12\cdot2^{-6}$.}\label{Rotating_Condensate_Table}
\end{center}
\end{table}

\subsection{Mott insulator to superfluid}
In our final experiment we use the RE-FEM, which was best-performing in our experiments, to present an interesting 2D simplification of the celebrated experiment by Greiner et al. \cite{Nature}. Before sketching the setup of the physical experiment, we stress that the computational complexity of the nonlinear methods was too high to carry out the computations at the desired resolution. The following numerical experiment therefore differs from the previous experiments in the sense that we will not compare the various methods and we do not have a numerical reference solution. Instead, we evaluate the numerical results by comparing them with the observations from the physical experiment. As we will see, the energy-conservative RE-FEM successfully captures the expected phase transition.

In the experiment by Greiner et al. \cite{Nature} an optical lattice is used to study quantum phase transitions of a gas of ultra cold bosons from the Mott insulator phase into the superfluid phase. In the Mott insulator phase, individual atoms are trapped in the sites of an optical lattice and the energy of particle interactions dominates over the kinetic energy. When the optical lattice is removed, the particles are free to move and the wave packages start to interfere with each other (phase transition). The aspect of physical interest is to study the subsequent interference pattern that is formed (or not formed) after the bosons are released from the lattice. As observed in \cite{Nature} the pattern is closely related to the strength of the lattice potential. More precisely, as the strength of the optical lattice is increased, the interference pattern after release from the lattice goes from a regular high-contrast pattern to an incoherent background. 
 
When the optical lattice potential is increased, the particle to particle interactions grow, at a certain threshold these interactions become too large for the gas to be described by a macroscopic wave function as in the Gross--Pitaevskii model. Instead, the (computationally very heavy) Bose--Hubbard model is needed for an accurate description of Mott insulators. In the following we shall assume that the energy of the particle interactions is high, but still below the threshold, so that we can use the Gross--Pitaevskii equation to compute an approximation of a Mott insulator ground state. In any case, after switching off the lattice, the Gross--Pitaevskii equation is the valid model for studying the dynamics of the ultra cold bosons and hence the formation of interference. A simulation faithful to the experiment by Greiner et al. in 3D could well require, of the order, $10^{12}$ degrees of freedom for the spatial discretization. As this is beyond computability, we will present a simplified 2D case with interesting dynamics requiring reasonable computational power.
 
From the physical setup described in \cite{Nature}, we can extract a configuration and a parameter range in which we expect to observe the phase transition from Mott insulator to superfluid together with a clear and structured interference pattern. As an admissible configuration for our numerical experiment, we set the strength of the
particle interactions to $\beta = 1000$ and select the following harmonic and optical lattice potentials:
\[V_{\text{h}}(x_1,x_2) = 2\hspace{2pt}(x_1^2+x_2^2),\ \  \ V_{\text{o}}(x_1,x_2) = 2000\hspace{2pt}(\sin(2\pi x_1)^2+\sin(2\pi x_2)^2) \]
The initial state (as an approximation for the Mott insulator) is taken to be the ground state associated with both the harmonic and the lattice potential, i.e.  $\uzero\in H^1_0(\Omega) $ solves  $-\frac{1}{2}\Delta u + (V_{\text{h}}+V_{\text{o}}+\beta|\uzero|^2)\uzero = \lambda \hspace{2pt} \uzero$, where $\lambda$ denotes the smallest eigenvalue. We solve the time-dependent Gross-Pitaevskii equation 
$\ci \partial_t u = -\frac{1}{2}\Delta u + V u + \beta |u|^2 u$ on $[-20,20]^2$ with homogeneous Dirichlet boundary conditions keeping only the harmonic potential, thus $V=2(x_1^2+x_2^2)$. We solve until final time $T=0.515$ with $h=0.01$, accounting for $16\cdot 10^6$ spatial degrees of freedom, and the number of time-steps is 2048, i.e. $\tau = 0.515/2^{11}$.

\begin{figure}[h!]
    \centering
    \begin{subfigure}[b]{0.45\textwidth}
        \includegraphics[width=1.1\textwidth]{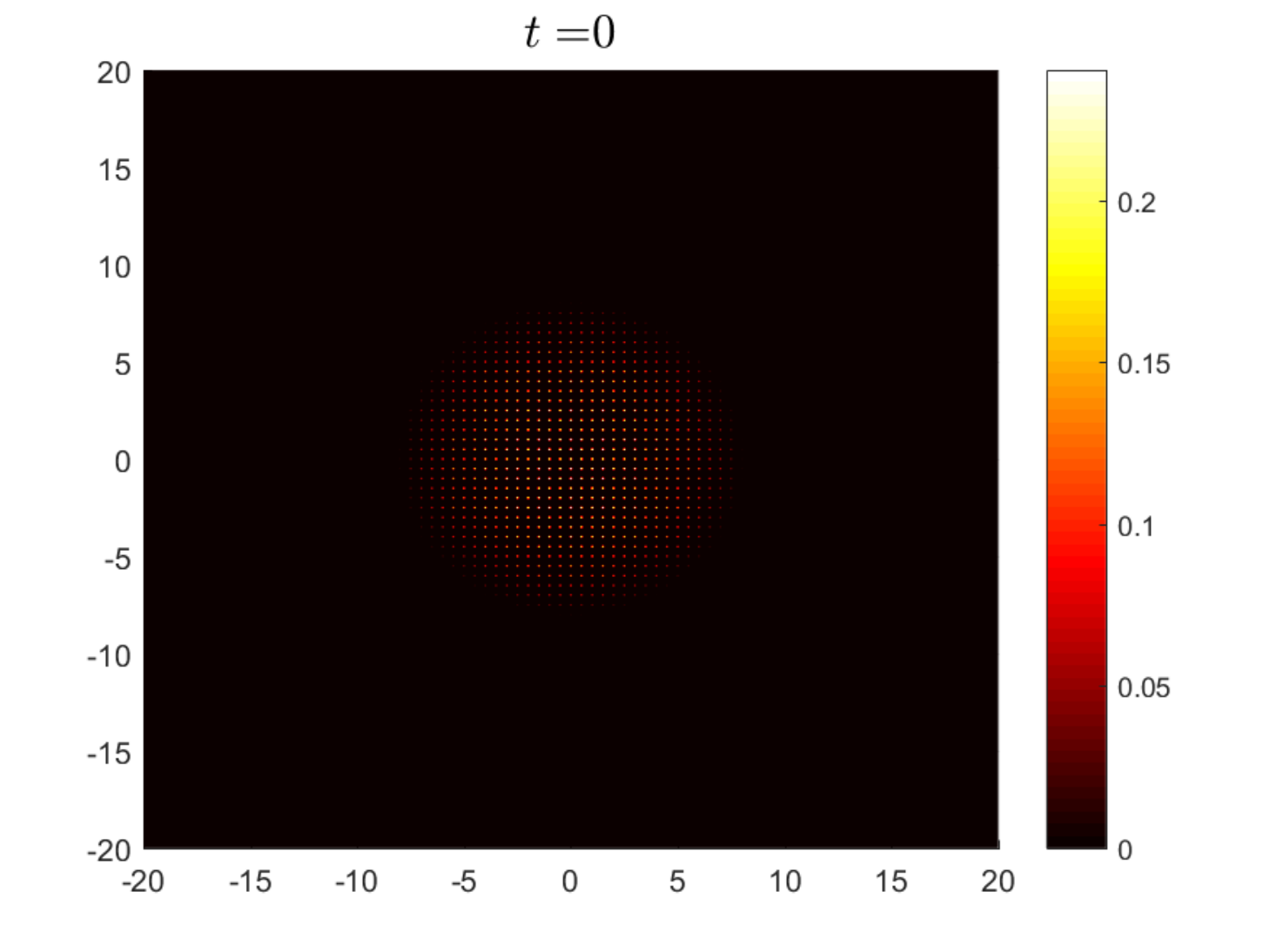}
        \caption{\scriptsize Ground state }
        \label{GroundState}
    \end{subfigure}
    ~  
    \begin{subfigure}[b]{0.45\textwidth}
        \includegraphics[width=1.1\textwidth]{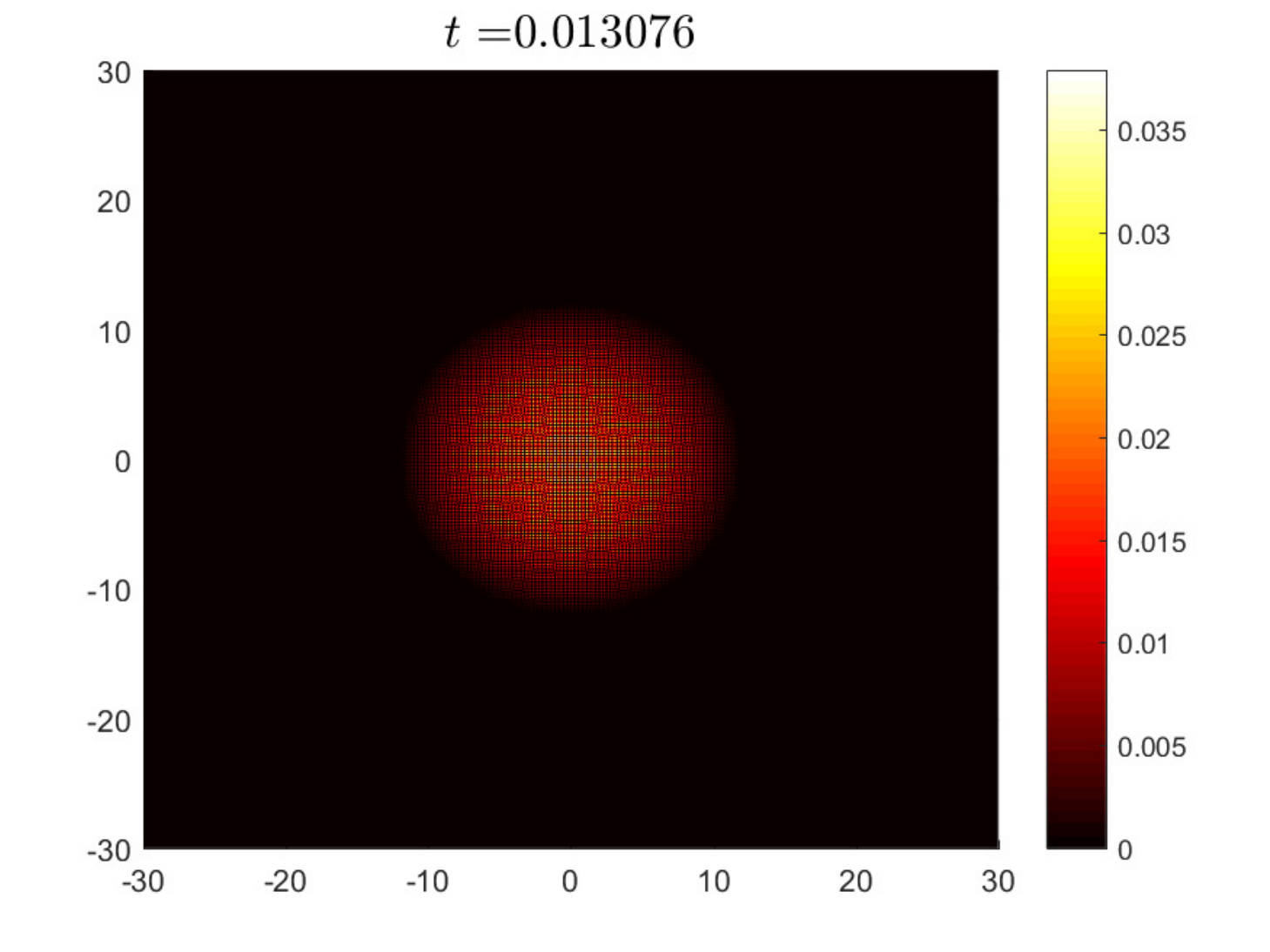}
        \caption{\scriptsize Early snapshot  }
        \label{EarlyTime}
    \end{subfigure}
    ~ 
~

    \begin{subfigure}[b]{0.45\textwidth}
        \includegraphics[width=1.1\textwidth]{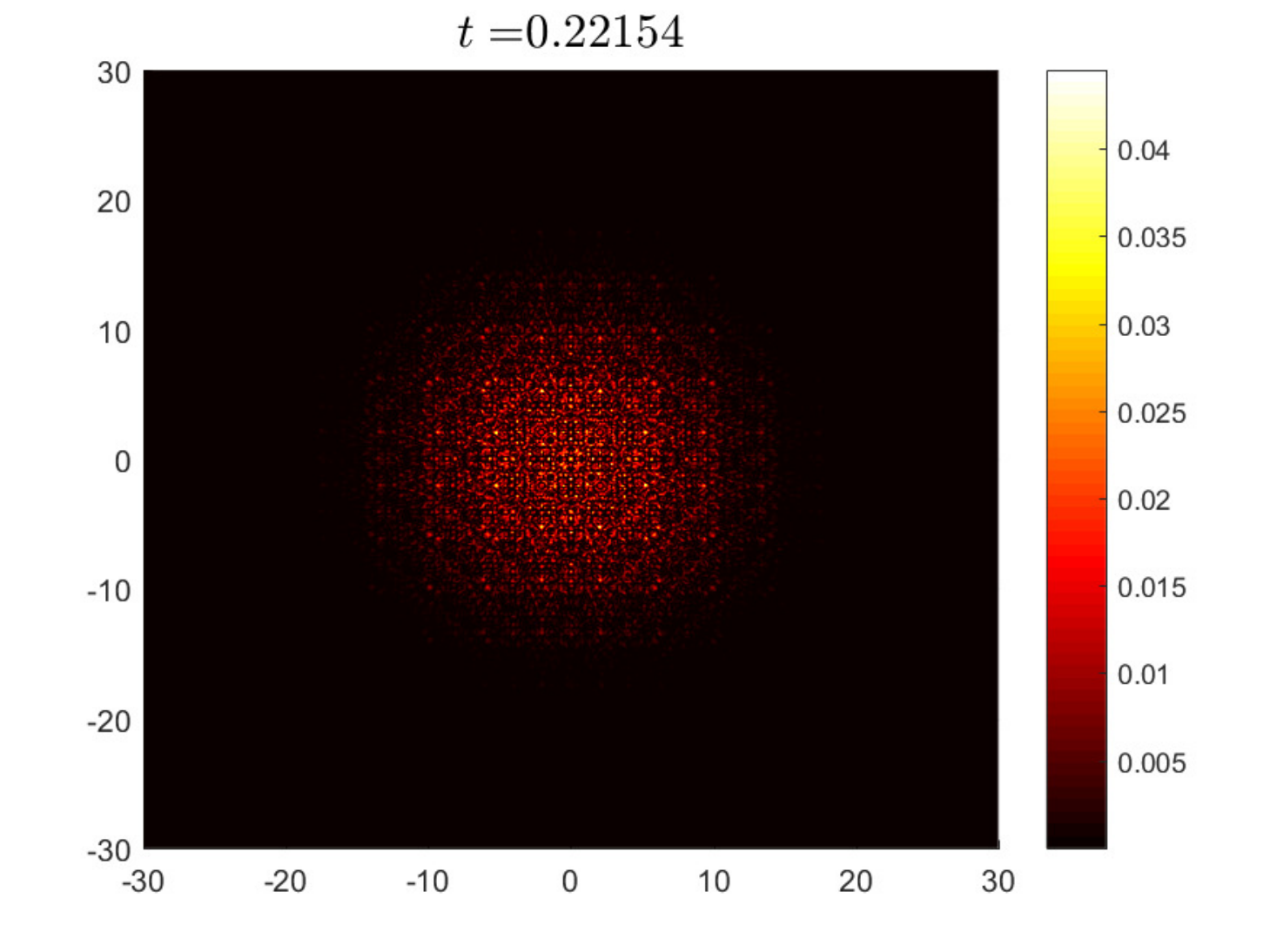}
        \caption{\scriptsize Half-time snapshot}
        \label{HalfTime}
    \end{subfigure}
    ~ 
    \begin{subfigure}[b]{0.45\textwidth}
        \includegraphics[width=1.1\textwidth]{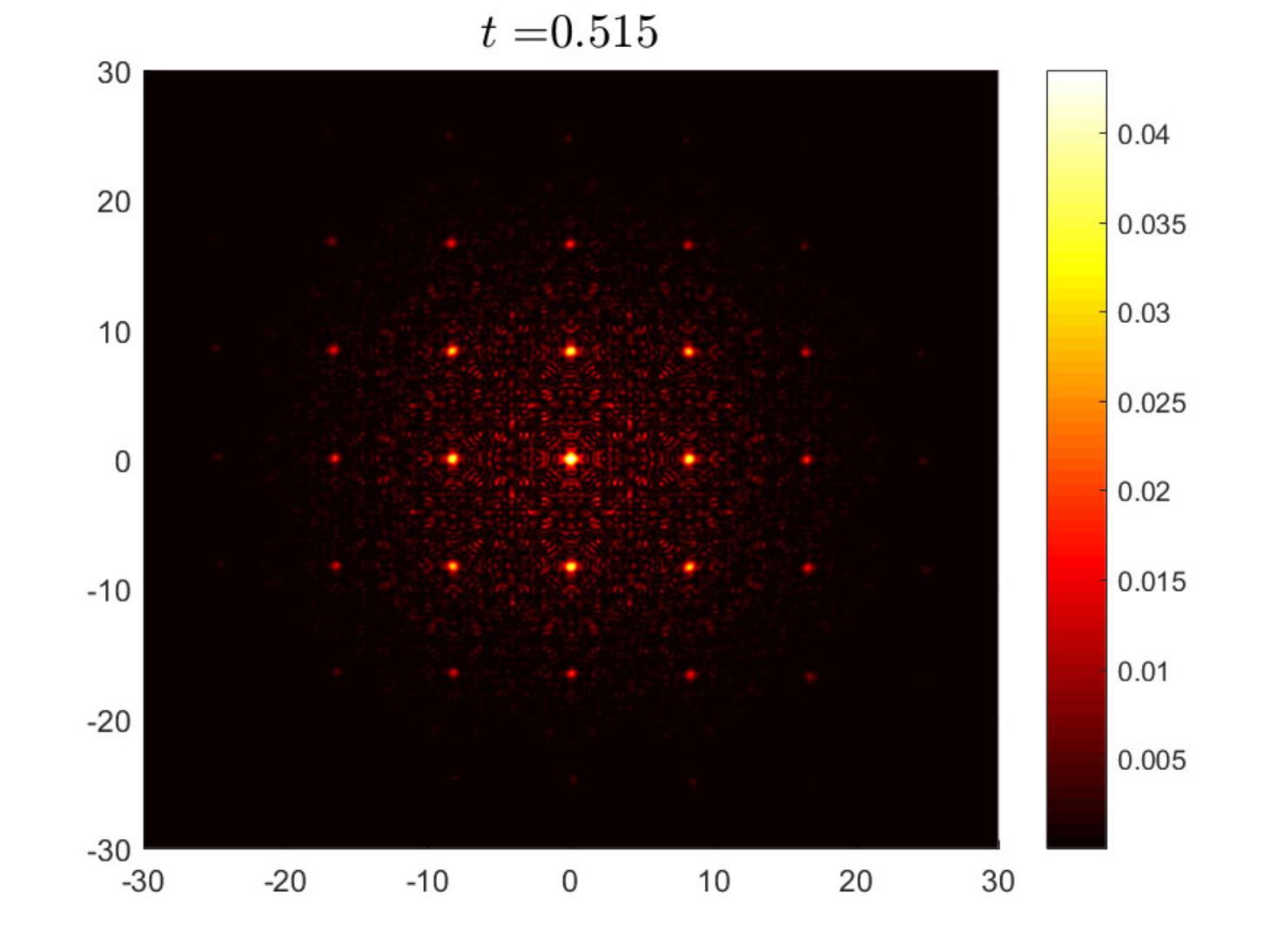}
        \caption{\scriptsize Solution at final time. }
        \label{FinalTime}
    \end{subfigure}
    ~ 
    \caption{\scriptsize Density plots, $|u_{h,\tau}|^2$, for $\tau = 0.515\cdot2^{-11}$ and $h=0.01$.  }
\end{figure}
 
As illustrated in Fig. \ref{GroundState}, the condensate starts in a regular pattern of localized and smooth clouds. After release, it quickly disintegrates into a more chaotic interference pattern with sharp spikes, illustrated in Fig. \ref{EarlyTime}. At $t=0.515$ highly localized and smooth clouds reappear in a regular lattice, shown in Fig. \ref{FinalTime}. The interference pattern shown in Fig. \ref{FinalTime} has precisely the structure of the pattern observed in the physical experiment by Greiner et al. We stress that this is only possible if the interference of the individual wave packages is accurately captured by the numerical scheme over many time steps. As we can see, the energy-conserving RE-FEM is able to reproduce the desired behavior.

\section{Conclusion} \label{Conclusions}

In this paper, we compared five different mass conservative time discretizations for the Gross-Pitaevskii equation. Two of the schemes (CN-FEM and RE-FEM) additionally preserve the energy, one scheme (IM-FEM) additionally preserves the symplectic structure and two schemes (LCN-FEM and Two-Step FEM) only preserve the mass. We observed that all schemes perform very well for simple experimental setups, where energy-conservation and symplecticity do not seem to play a crucial role to obtain good results. However, when moving towards more complicated setups with reduced regularity and additional potential terms, the conservation of energy becomes more and more important and the purely mass-conservative methods, LCN-FEM and Two-Step FEM, can fail to produce accurate approximations on reasonable scales. We also observed that in our experiments the energy-conserving methods gave visibly better approximations than the symplectic IM-FEM.

None of our experiments gave any indications that convergence could depend on a possible coupling between mesh size and time-step size as often seen in the literature to prove convergence. Still we saw that an unlucky choice of the ratio $h/\tau$ can cause an energy blow-up for the LCN-FEM and the Two-Step FEM. This blow-up vanishes slowly for $h,\tau \rightarrow 0$, hence not contradicting convergence, however it can introduce severe practical constraints. Surprisingly, we also find that the symplectic IM-FEM can produce completely erroneous density plots, whereas the energy-conservative methods do not. In short, the energy-conservative methods appear to be more reliable in complex physical settings. The linear energy-conservative RE-FEM achieves errors at least on par with the nonlinear CN-FEM, thereby making it by far the most computational efficient method.

\newpage
\providecommand{\href}[2]{#2}
\providecommand{\arxiv}[1]{\href{http://arxiv.org/abs/#1}{arXiv:#1}}
\providecommand{\url}[1]{\texttt{#1}}
\providecommand{\urlprefix}{URL }

\end{document}